\documentclass[11pt]{article}
\usepackage{epsfig}
\usepackage{amssymb}
\usepackage{amsmath}
\usepackage{array}
\def\BbbR{{\Bbb R}}
\def\BbbN{{\Bbb N}}   
\def\BbbZ{{\Bbb Z}}
\def\BbbC{{\Bbb C}}
\def\BbbR{{\Bbb R}}

\def\BbbH{{\Bbb H}}

\newtheorem{theo}{Theorem}
\newtheorem{lemm}{Lemma}
\newtheorem{coro}{Corollary}

\newtheorem{prop}{Proposition}
\begin{document}
\title{Exploring Surfaces through Methods from the Theory of Integrable Systems.
Lectures on the Bonnet Problem\footnote{Lectures given at the School on Differential Geometry on 12-30 April 1999 at the
Abdus Salam International Centre for Theoretical Physics (ICTP), Trieste}}
\author{Alexander I. Bobenko
\\ Fachbereich Mathematik, Technische Universit\"{a}t
Berlin,
\\ Strasse des 17. Juni 136, 10623 Berlin, Germany}
\date{}
\maketitle
\begin{abstract}
 A generic surface in Euclidean 3-space is determined uniquely by its metric and 
curvature. Classification of all special surfaces where this is not the case, i.e. of surfaces
possessing isometries which preserve the mean curvature, is known as the Bonnet problem. 
Regarding the Bonnet problem, we show how analytic methods of the theory of integrable 
systems -- such as finite-gap integration, isomonodromic deformation,  and loop 
group description --
can be applied for studying global properties of special surfaces. 
\end{abstract}

\newpage 
\section{Lecture. Quaternionic Description of Surfaces. 
Bonnet Problem}							\label{s1}

\subsection{Differential equations of surfaces}			\label{s1.1}
Let ${\cal F}$ be a smooth orientable surface in 3-dimensional Euclidean space.
The Euclidean metric induces a metric $\Omega$ on this surface, which 
in turn generates the complex structure of a Riemann surface $\cal R$. 
 Under 
such a pa\-ra\-met\-ri\-za\-tion, which is called {\em conformal}, the surface
$\cal F$ is given by an immersion
\[F=(F_1, F_2, F_3) : {\cal R} \rightarrow \BbbR^3,\]
and the metric is conformal: $\Omega = e^{u}\, dzd\bar{z}$, where $z$ is
a local coordinate on $\cal R$. 

One should keep in mind that a complex coordinate is defined up to holomorphic $z\to w(z)$  
transformation. This freedom will be used to simplify the corresponding
equations.

The conformal pa\-ra\-met\-ri\-za\-tion gives the following normalization of 
$F(z, \bar{z})$:
\begin{eqnarray}
<F_z, F_z>=<F_{\bar{z}}, F_{\bar{z}}>=0,\  
<F_z, F_{\bar{z}}>=\frac{1}{2} e^{u},                          \label{R31.2.1}
\end{eqnarray}
where the brackets denote the scalar product
\[<a,b>=a_1 b_1+ a_2 b_2+a_3 b_3,\]
and $F_z$ and $F_{\bar{z}}$ are the partial derivatives $\frac{\partial F}{\partial
z}$ and $\frac{\partial F}{\partial\bar{z}}$, where
$$
 \frac{\partial}{\partial z}=\frac{1}{2}
\left(\frac{\partial}{\partial x}-i\frac{\partial}{\partial y}\right),\qquad
\frac{\partial}{\partial \bar{z}}=\frac{1}{2}
\left(\frac{\partial}{\partial x}+i\frac{\partial}{\partial y}\right). 
$$
The vectors $F_z, F_{\bar{z}}$, as well as the normal $N$, with
\begin{eqnarray}
<F_z,N>=<F_{\bar{z}},N>=0, \ \ <N,N>=1,                        \label{R31.2.2}
\end{eqnarray}
define a moving frame on the surface, which due to (\ref{R31.2.1}, 
\ref{R31.2.2}) satisfies
the  following Gauss-Weingarten equations:
\begin{eqnarray}
&\sigma_z={\cal U}\sigma,\ \sigma_{\bar{z}}={\cal V}\sigma,\ 
\sigma=(F_z,F_{\bar{z}},N)^T,
                                                  \label{R31.2.3}\\
&						  \nonumber\\
&\cal U=\left(\begin{array}{ccc}
u_z & 0 & Q\\
0 & 0 & \frac{1}{2}He^u\\
-H & -2e^{-u}Q & 0
\end{array}\right),\                          
\cal V=\left(\begin{array}{ccc}
0 & 0 &  \frac{1}{2}He^u\\
0 & u_{\bar{z}} &  \bar{Q}\\
-2e^{-u}\bar{Q} & -H & 0
\end{array}\right),                                     \label{R31.2.4}           \end{eqnarray}
where
\begin{eqnarray}
Q=<F_{zz},N>, \ \ <F_{z\bar{z}},N>=\frac{1}{2}He^{u}.    \label{R31.2.5}
\end{eqnarray}
The quadratic differential $Q dz^2$ is called the {\em Hopf differential}.
The first and the second fundamental forms 
\begin{eqnarray*}					
< dF, dF >&=&
< I  { dx \choose dy},\
{ dx \choose dy } >,\ z=x+iy,\\
-< dF, dN >&=&
< II{ dx \choose dy },
{ dx \choose dy } >
\end{eqnarray*}
are given by the matrices
\begin{eqnarray}					\label{fund.forms}
I=e^u\left(\begin{array}{cc}1&0\\0&1 \end{array}\right),\ 
II=\left(\begin{array}{cc} Q+\bar{Q}+He^u&i(Q-\bar{Q})\\
i(Q-\bar{Q})&-(Q+\bar{Q})+He^u \end{array}\right).
\end{eqnarray}
The {\em principal curvatures} $k_1$ and $k_2$ are the eigenvalues of the matrix
\hbox{$II \cdot I^{-1}$}. This gives the following expressions for the {\em mean} and 
the {\em Gaussian curvatures}:
\begin{eqnarray*}
&H=\frac{1}{2}\ (k_1+k_2)=\frac{1}{2}\ {\rm tr}\ (II\cdot I^{-1}),\\
&K=k_1 k_2=\det\ (II\cdot I^{-1})=H^2-4Q\bar{Q}e^{-2u}.
\end{eqnarray*}

A point $P$ of the surface $\cal F$ is called {\em umbilic} if the principal
curvatures at this point coincide $k_1(P)=k_2(P)$. 
The Hopf differential vanishes $Q(P)=0$ exactly at umbilic points of the
surface. 

Coordinates  in which both fundamental forms are diagonal are
called {\em curvature line coordinates} and the corresponding parametrization 
(not necessarily conformal) is
called a {\em curvature line parametrization}. A curvature line
parametrization always exists in a neighborhood of a non-umbilic point. 
Near umbilic points, curvature lines form more complicated patterns. 

The Gauss--Codazzi equations, which are the compatibility conditions of 
equations (\ref{R31.2.3}, \ref{R31.2.4}),
\[{\cal U}_{\bar{z}}-{\cal V}_z+[{\cal U},{\cal V}]=0, \]
have the following form:
\begin{equation}		\label{complex version of Gauss-Codazzi}
\begin{array}{lrcl}	
 {\rm Gauss\  equation}\qquad &                    
	u_{z \bar z} + \frac{1}{2}\,
	H^2\,e^u - 2 \vert Q\vert^2\,e^{-u} 
	& = & 0,
\\
 {\rm Codazzi\  equation}\qquad &
	Q_{\bar z} & = & \frac{1}{2}\,H_z\,e^u.
\end{array}
\end{equation}

These equations are necessary and sufficient for existence of the corresponding
surface.

\begin{theo} {\bf (Bonnet theorem)}.				\label{Bonnet.t}
Given a metric $e^u\, dzd\bar{z}$, a quadratic differential $Q\, dz^2$, and a
function $H$ on $\cal R$ satisfying the Gauss--Codazzi equations, there exists an
immersion
$$
F:\tilde{\cal R} \to \BbbR^3
$$ 
with the fundamental forms (\ref{fund.forms}). Here $\tilde{\cal R}$ is the
universal covering of $\cal R$. The immersion $F$ is unique up to Euclidean
motions in $\BbbR^3$.
\end{theo}

We finish this section with some basic facts about a special class of 
surfaces.
A conformal curvature line parametrization is called {\em isothermic}.
In this case the preimages of the curvature lines are the lines
$x={\rm const}$ and $y={\rm const}$ on the parameter domain, where
$z=x+iy$ is a conformal coordinate.
A surface  is called {\em isothermic}
if it allows isothermic parametrization. Isothermic surfaces are 
divided by their curvature lines into infinitesimal squares.
Written in terms of an isothermic coordinate $z$ the Hopf differential 
$Q(z, \bar{z}) dz^2$ of an isothermic surface is real, 
i.e. $Q(z, \bar{z})\in \BbbR$.

In terms of arbitrary conformal coordinates, isothermic surfaces can be 
characterized as follows.
 
\begin{lemm}\label{Remark: how to indicate isothermic}
Let $F:{\cal R}\to {\BbbR}^3$ be a conformal immersion of an umbilic free surface in 
${\BbbR}^3$. The surface is isothermic if and only if 
there exists a holomorphic non-vanishing differential
$f(z) dz^2$ on $\cal R$ and a function 
$q:{\cal R}\to\BbbR_*$ such that the Hopf differential 
is of the form
\begin{equation}		\label{eq: how to indicate isothermic}
	Q(z,\bar z) = f(z)\,q(z,\bar z).
\end{equation}
\end{lemm}
It is easy to see that $w=\int \sqrt{f(z)}\,d z$ is an isothermic 
coordinate.

\subsection{Quaternionic description of surfaces}		\label{s1.2}

We construct and investigate surfaces in ${\BbbR}^3$ by analytic methods. 
For this purpose it is convenient to use the Lie algebra isomorphism 
$so(3)=su(2)$ and to rewrite the equations 
(\ref{R31.2.3}, \ref{R31.2.4}) for the moving frame in terms of 2 by 2 matrices.
This quaternionic description turns out to be 
useful for analytic studies of general curves and surfaces in 3- and 
4-spaces as well as for investigation of special classes of surfaces 
\cite{BoCMC, KuS2, DPW, Bo2x2, KPP, PP}.

Let us denote the algebra of quaternions by $\BbbH$, the multiplicative quaternion
group by ${\BbbH_*}={\BbbH}\setminus\{0\}$, and their standard basis by $\{{\bf1}
,{\bf i},{\bf j},{\bf k}\}$, where
\begin{eqnarray}
{\bf i}{\bf j}={\bf k},\ 
{\bf j}{\bf k}={\bf i},\
{\bf k}{\bf i}={\bf j}.                                         \label{R32.2.9}
\end{eqnarray}
This basis can be represented by the Pauli matrices $\sigma_\alpha$ as follows:
\begin{eqnarray} 
&\sigma_1=\left( \begin {array}{cc} 0&1\\1&0 \end{array}\right)=i\ {\bf i},\ \
\sigma_2=\left( \begin {array}{cc} 0&{-i}\\i&0 \end{array}\right)=i\ {\bf j},\nonumber\\
&\sigma_3=\left( \begin {array}{cc} 1&0\\0&{-1} \end{array}\right)=i\ {\bf k},\ \
{\bf1}=\left( \begin {array}{cc} 1&0\\0&1 \end{array}\right).     \label{R32.2.10}
 \end{eqnarray}
We identify
$\BbbH$ with 4-dimensional Euclidean space
$$
q=q_0 {\bf 1}+ q_1 {\bf i}+q_2 {\bf j}+q_3 {\bf k}
\longleftrightarrow \ q=(q_0, q_1, q_2, q_3)\in\BbbR ^4.
$$
The sphere $S^3\subset\BbbR^4$ is then naturally identified with the group of unitary
quaternions $\BbbH_1=SU(2)$.
3-dimensional Euclidean space is identified with the space of imaginary quaternions 
${\rm Im}\ {\BbbH}$
\begin{eqnarray} 
X =-i\sum_{\alpha =1}^3 X_\alpha \sigma_\alpha \in{\rm Im}\ {\BbbH}\ 
\longleftrightarrow \ X=(X_1, X_2, X_3)\in\BbbR ^3.                             \label{R32.2.11}
\end{eqnarray}
The scalar product of vectors in terms of quaternions and matrices is then
\begin{eqnarray} <X,Y>=-\frac{1}{2}(XY+YX)=-\frac{1}{2}{\rm tr}\, XY.            \label{R32.2.12}
\end{eqnarray}
We will also denote by $F$ and $N$ the matrices obtained in this way from the 
vectors $F$ and $N$.

Let us take $\Phi \in \BbbH_*$ which transforms the basis 
${\bf i}, {\bf j}, {\bf k}$  into the frame $F_x, F_y, N$:
\begin{equation}
F_x=e^{u/2}\Phi^{-1}{\bf i}\Phi,\
F_y=e^{u/2}\Phi^{-1}{\bf j}\Phi,\
N=\Phi^{-1}{\bf k}\Phi.                                         \label{R32.2.14}
\end{equation}
Then 
\begin {eqnarray}
F_z=-ie^{u/2}\Phi^{-1}\left(\begin{array}{cc}
0 & 0 \\
1 & 0 \end{array}\right)\Phi,\  \
F_{\bar {z}}=-ie^{u/2}\Phi^{-1}\left(\begin{array}{cc}
0 & 1 \\
0 & 0 \end{array}\right)\Phi,                                    \label{R32.2.15} 
\end {eqnarray}
and all the conditions (\ref{R31.2.1}) are automatically satisfied.

The quaternion $\Phi$ satisfies linear differential equations. To derive them we
introduce matrices
\begin {equation}
U=\Phi_z \Phi^{-1},\qquad V=\Phi_{\bar z} \Phi^{-1}.                 \label{R32.2.16}
\end {equation}
The compatibility condition $ F_{z\bar z}=F_{\bar z z}$ for (\ref{R32.2.15}) implies
\begin {eqnarray*}
V_{22}-V_{11}={u_{\bar z}\over 2},\ \ U_{11}-U_{22}={u_z \over 2},\ \ 
U_{21}=-V_{12},
\end {eqnarray*}
where $U_{kl}$ and $V_{kl}$ are the matrix elements of $U$ and $V$. In the same
way one obtains from (\ref{R32.2.15}) 
\begin {eqnarray*}
F_{z {\bar z}}=\frac{1}{2} H e^uN\ &\rightarrow & \  
U_{21}=-V_{12}=\frac{1}{2}H e^{u/2} \\
F_{zz}=u_z F_z+ Q N\ &\rightarrow & \
U_{12}=-Qe^{-u/2}\\
F_{\bar z \bar z}=u_{\bar z} F_{\bar z}+ {\bar Q} N\ &\rightarrow & \
V_{21}=\bar Qe^{-u/2}.
\end {eqnarray*}
Recall that $\Phi$ is defined up to multiplication by a scalar factor. We normalize
this factor by the condition
\begin{equation}						\label{SpinNorm}
\det \Phi =e^{u\over 2},
\end{equation}
for reasons which will be clarified in the next section. For the traces of $U$ and $V$
this implies
$$
{\rm tr}\, U={u_z\over 2},\qquad {\rm tr}\, V={u_{\bar{z}}\over 2}.
$$

Finally one arrives at the following

\begin{theo} 							\label{1.t2}
Using the isomorphism (\ref{R32.2.11}), the moving frame $F_z, F_{\bar z}, N$
of a conformally para\-met\-rized surface ($z$ is a conformal coordinate) is described
by formulas (\ref{R32.2.14}),(\ref{R32.2.15}), where $\Phi\in \BbbH_*$ 
satisfies the equations
(\ref{R32.2.16}) with $U$, $V$ of the form
\begin {eqnarray}
U=\left( \begin {array}{cc}\displaystyle\frac {u_z}{2} & -Q e^{-u/2}\\
\displaystyle \frac {1}{2} H e^{u/2}  &  0  \end{array}
 \right),\ \  
V=\left( \begin {array}{cc} 0  &\displaystyle  -\frac {1}{2} H e^{u/2}\\
\bar Q e^{-u/2}  & \displaystyle  \frac {u_{\bar z}}{2}  \end{array}\right).   
                                                                    \label{R32.2.18}
\end {eqnarray}
\end{theo}

\begin{coro}							\label{1.c1}
The conformal frame $\Phi$ satisfies the Dirac equation
\begin{equation}						\label{Dirac1}
e^{-u/2}
\left( \begin {array}{cc} 0 & \partial_z\\
-\partial_{\bar{z}}  &  0  \end{array}
 \right)\Phi ={1\over 2}H\Phi.
\end{equation} 
\end{coro}

It turnes out that at this point the whole construction can be reversed.
Namely, starting with a solution to the Dirac equation one can derive a Weierstrass
type representation (see (\ref{Weier}) below) for conformally parametrized surfaces.
This idea was recently developed by Konopelchenko \cite{Kon} and further in 
\cite{Tai, PP, KuS2}, although in other forms the Weierstrass representation of
surfaces was known already to Eisenhart \cite{Eis} and Kenmotsu \cite{Ken}.

\begin{theo}							\label{1.t3}
Let $D\subset \BbbC$ be a simply connected domain and  
$(s_1, \bar{s}_2)^T :D\to\BbbC^2$ be a solution to the Dirac equation
with the potential $p\in C^\infty(D)$
\begin{equation}						\label{Dirac2}
\left( \begin {array}{cc} 0 & \partial_z\\
-\partial_{\bar{z}}  &  0  \end{array}
 \right)
 \left( \begin {array}{c} s_1\\ \bar{s}_2  \end{array}
 \right) = p\left( \begin {array}{c} s_1\\ \bar{s}_2  \end{array}
 \right).
\end{equation}
Then
\begin{equation}						\label{Phi}
\Phi=
\left( \begin {array}{cc} s_1 & -s_2\\
\bar{s}_2  &  \bar{s}_1  \end{array}
 \right): D\to \BbbH_*
\end{equation}
is a conformal frame (\ref{R32.2.15}) of the conformally immersed surface
\begin{eqnarray}
F_1+iF_2 &=& \int s_1^2\, dz - \bar{s}_2^2\, d\bar{z}      \nonumber\\
F_3 &=& \int s_1 s_2\, dz+ \bar{s}_1 \bar{s}_2\, d\bar{z}. 	\label{Weier}
\end{eqnarray}
The metric and the mean curvature of the surface are given by
\begin{equation}
e^u\, dzd\bar{z}=(\mid s_1\mid ^2+\mid s_2\mid ^2)^2\, dzd\bar{z},\quad
H=2 p e^{-u/2}.							\label{WeierFF}
\end{equation}
\end{theo}

{\em Proof.} Note that $(-s_2, \bar{s}_1)^T$ is also a solution to (\ref{Dirac2})
due to the symmetry of the Dirac equation. At this point $\Phi$ given by (\ref{Phi}) 
can be identified with the conformal frame $\Phi$ of Corollary \ref{1.c1}. 
The formula for the metric (\ref{WeierFF}) follows from (\ref{SpinNorm}).
Substituting it into our previous formulas (\ref{R32.2.15}) for conformal frame one
defines
$$
F_z:=-i\left( \begin {array}{cc} s_1 s_2 & -s_2^2\\
s_1^2  &  -s_1 s_2  \end{array}
 \right),\quad
F_{\bar{z}}:=-i\left( \begin {array}{cc} \bar{s}_2 \bar{s}_1 & \bar{s}_1^2\\
-\bar{s}_2^2  &  -\bar{s}_2\bar{s}_1  \end{array}
\right).
$$ 
These formulas are automatically compatible. Integrating them one arrives at (\ref{Weier}).
 
\subsection{Spinor description of surfaces}			\label{s1.4}

As shown in \cite{Bo2x2}, the quaternionic description of the previous section is 
actually a global one.
Let $\cup_i D_i={\cal R}$ be an open covering of $\cal R$ with local coordinates
$z_i:D_i\to \BbbC$. Conditions (\ref{R32.2.15}, \ref{SpinNorm}) determine a quaternionic
valued smooth $\Phi(z_i,\bar{z}_i)$ uniquely up to sign on each $D_i$. 
To establish the global nature of $\Phi$ recall that a holomorphic 
line bundle $S$ is called a {\em spin
bundle} if it satisfies $S\otimes S=K$, where $K$ is the canonical bundle. 

Denote the first column of $\Phi$ by 
$$
\left( \begin {array}{c} S_1 \\ \bar{S}_2  \end{array} \right).
$$

\begin{lemm}
$S_1$ and $S_2$ are smooth sections of the same holomorphic spin bundle $S$.
\end{lemm}

{\em Proof.} Consider two intersecting $D_i\cap D_j\not= \emptyset$ with corresponding 
$\Phi_i(z_i,\bar{z}_i)$ and $\Phi_j(z_j,\bar{z}_j)$. Identifying the representations for 
the Gauss map in terms of $\Phi_i$ and $\Phi_j$ one obtains on $D_i\cap D_j$
$$
\Phi_i=
\left( \begin {array}{cc} c_{ij} & 0\\
0  &  \bar{c}_{ij}  \end{array}
 \right) \Phi_j.
$$
with some $c_{ij}:D_i\cap D_j\to \BbbC_*$. Further, identifying the tangent frames 
$$
F_{z_i}=F_{z_j}{dz_j\over dz_i}
$$ 
and using 
$\Phi\sigma_2\Phi^T\sigma_2=\det\Phi$ one obtains
$$
\Phi_i^T
\left( \begin {array}{cc} 1 & 0\\
0  &  0  \end{array} \right)\Phi_i={dz_j\over dz_i}
\Phi_j^T
\left( \begin {array}{cc} 1 & 0\\
0  &  0  \end{array} \right)\Phi_j,
$$
which finally implies
$$
c_{ij}^2={dz_j\over dz_i}.
$$
The transition functions $c_{ij}:D_i\cap D_j\to \BbbC_*$ defined through $\Phi_i$ obviously 
satisfy the cocycle condition $c_{ij}c_{jk}=c_{ik}$ and thus define a line bundle $S$ with 
$S\otimes S=K$.

In local coordinates $S_n$ are written naturally as $S_n=s_n(z_i,\bar{z}_i)\, \sqrt{dz_i}$. 

Using the  equivalence of spinor
representation of conformal frames of surfaces and solutions of the 
Dirac equation, 
proven in Corollary \ref{1.c1} and Theorem \ref{1.t3}, we
arrive at the following global reformulation \cite{Tai, PP} of Theorem \ref{1.t3}.

\begin{theo}							\label{1.t4}
 A half-density $p$ (i.e. a smooth section of 
$K^{1\over 2}\otimes \bar{K}^{1\over 2}$) and two not simultaneously vanishing 
spinors $S_1, S_2$  (i.e. smooth sections of $S\cong K^{1\over 2}$ with
$(S_1, S_2)\not= (0, 0) \ \forall P\in{\cal R}$) satisfying the Dirac equation 
(\ref{Dirac2}) determine through
\begin{eqnarray}
F_1+iF_2 &=& \int S_1^2 - \bar{S}_2^2      \nonumber\\
F_3 &=& \int S_1 S_2+ \bar{S}_1 \bar{S}_2 			\label{weierGlobal}	
\end{eqnarray}
a conformal immersion $F:{\tilde{\cal R}\to \BbbR^3}$, where $\tilde{\cal R}$ is the universal
covering of $\cal R$. The metric and the mean curvature of the immersion are given by
$$
e^u\, dzd\bar{z}=(\mid S_1\mid ^2+\mid S_2\mid ^2)^2,\quad
H=2 p e^{-u/2}.							
$$
\end{theo}

{\it Remark}. In the case of minimal surfaces $H=0$ the spinors $S_1$ and $S_2$ are
holomorphic and the representation (\ref{weierGlobal})
is known as the spinor Weierstrass representation \cite{Sul, Bo2x2, KuS1}.
\vspace{0.5cm}

On a Riemann surface of genus $g$ there exist $2^{2g}$ non-isomorphic spin bundles which are
distinguished by different {\em spin structures}. For a geometric interpretation of the spin
structure of the spin bundle $S$ in terms of the immersion (\ref{weierGlobal}) we refer to 
\cite{Bo2x2}. Spin structures classify regular homotopies of immersions \cite{Pin}.

\subsection{Alternative descriptions of surfaces and the Bonnet problem}					\label{s1.5}			
Bonnet Theorem \ref{Bonnet.t} characterizes surfaces via the coefficients 
$e^u, H, Q$ of their fundamental forms. These coefficients are not independent and 
are subject to the Gauss--Codazzi equations (\ref{complex version of Gauss-Codazzi}). 
A natural question is whether some of these data are superfluous. The following 
natural candidates for more "economic" characterization of surfaces were studied.
\vspace{0.5cm}

(i) The most geometric setting of the problem is the oldest one and is due to Bonnet.
He posed the question whether one can eliminate the Hopf differential from the description of
surfaces, i.e. whether the {\em metric} and the {\em mean curvature function} 
$$
e^u, \qquad H
$$
alone suffice to describe a surface completely. Generic surfaces are determined uniquely by the 
metric and the mean curvature function. Bonnet himself \cite{Bonnet} made the initial progress in 
the investigation of the special surfaces where it is not the case, i.e. which possess
non-congruent isometric "relatives" with the same mean curvature function.
The rest of these lectures is devoted to this problem, which is fairly named the 
{\em Bonnet problem}.
\vspace{0.5cm}

(ii) The {\em conformal Hopf differential} 
$$
q:= Q e^{-{u\over 2}}.
$$
Note that whereas the Hopf differential is a quadratic differential, i.e. a section of the line
bundle $K^2$, the conformal Hopf differential is more exotic - it is a section of 
$K^{3\over 2}\otimes \bar{K}^{-{1\over 2}}$.
The reason for its introduction  by U. Pinkall
is that $q$ is invariant with respect to the  M\"{o}bius 
transformations of the ambient $\BbbR^3$. A non-isothermic surface is uniquely 
determined by $q$ up to M\"{o}bius transformations. Counting dimensions, one immediately 
observes that generic sections of $K^{3\over 2}\otimes \bar{K}^{-{1\over 2}}$ do not correspond
to surfaces in $\BbbR^3$.  A proper equation for $q$ of surfaces in $\BbbR^3$ is still unknown.
\vspace{0.5cm}

(iii) The {\em Dirac potential} or {\em mean-curvature half-density} 
$$
p={1\over 2}He^{-{u\over 2}}.
$$ 
As one can see from its definition, this potential is a half-density, i.e. a section of the 
line bundle $K^{1\over 2}\otimes \bar{K}^{1\over 2}$. Recently, description 
of surfaces through Dirac spinors attracted much attention \cite{KuS2, Tai, PP}. Unfortunately, 
one has neither existence nor uniqueness in this description. A generic Dirac operator (with 
generic potential) has trivial kernel, thus generic half-densities do not yield surfaces. On 
the other hand, there may exist many immersions with the same potential $p$, for example all
special surfaces appearing in the Bonnet problem. 
 \vspace{0.5cm}

Returning to the Bonnet problem, note that already Bonnet indicated all special surfaces 
which possess non-congruent isometric "relatives" with the same mean curvature function.
There are three cases when this happens.
\vspace{0.5cm} 

{\em 1. Constant mean curvature surfaces}. Let ${\cal F}$ be a surface with
constant mean curvature $H$. The Gauss--Codazzi equations 
(\ref{complex version of Gauss-Codazzi}) are obviously invariant with respect to
the transformation
$$
Q \to Q_t = e^{it} Q,\qquad t\in\BbbR.
$$ 
Applying the Bonnet theorem one obtains the one parameter family 
${\cal F}_t,\ {\cal F}={\cal F}_0$ of
isometric surfaces with the same constant mean curvature $H$.
In the last ten years there was much interest in studying global properties of
surfaces with constant mean curvature and now they are rather well 
investigated by various methods (see for example \cite{Wen, PS, Kap, BoCMC, KGBKS}) 
including methods of the theory of integrable systems.  
\vspace{0.5cm}

{\em 2. Bonnet pairs} are exactly two non-congruent isometric surfaces
${\cal F}'$ and ${\cal F}''$ with the same mean curvature function. The theory
of Bonnet pairs is very closely related \cite{BianchiBonnet, KPP} to the 
theory of isothermic surfaces and as such belongs also to geometry described by
integrable systems. Up to now, global theory of Bonnet pairs is not well
developed, in particular it is unknown whether there exist compact Bonnet
pairs, a question first posed in \cite{LawTri}.
\vspace{0.5cm}

{\em 3. Bonnet families}. In \cite{Bonnet}, Bonnet himself was able to show that
besides the surfaces with constant mean curvature there exists a class of 
surfaces, depending on finitely many parameters which possess one-parameter 
family of isometries preserving the mean curvature. These surfaces
were studied by many authors \cite{Haz, Gra, Car, Che, BE1, Rou} and recently global 
classification \cite{BE2} of them was obtained using methods from the theory of 
integrable systems.
\vspace{0.5cm}

The remaining three lectures are devoted to consideration of these three cases.

\newpage 
\section{Lecture. Constant mean curvature surfaces}		\label{s2}
The content of this lecture is based on results obtained in \cite{BoCMC}.
\subsection{Associated family}					\label{s2.1}

If the mean curvature of $\cal F$ is constant, then the Gauss-Codazzi
equations
$$
u_{z\bar{z}} +\frac{1}{2}H^2e^u-2Q\bar{Q}e^{-u}=0,\qquad Q_{\bar{z}}=0,
$$
are invariant with respect to the transformation
\begin{eqnarray}
Q\to Q^{t}=\Lambda Q,\quad |\Lambda|=1,                           \label{CMC1.5.1}
\end{eqnarray}
Integrating the equations for the moving frame with
the  coefficient $Q$ replaced by $Q^t=\Lambda Q$ we obtain a one-parameter
family ${\cal F}^t$ of surfaces. All the 
surfaces ${\cal F}^t$ are isometric and have the same constant mean 
curvature. Treating $t$ as a deformation parameter we obtain the first family of
special surfaces indicated by Bonnet (see Section \ref{s1.5}).

\begin{theo}   																																										
\label{CMC1.t5.1} 
  Every constant mean curvature surface has a 
one-pa\-ra\-me\-ter family of isometric deformations preserving the
mean curvature. The deformation is described by 
 the transformations (\ref{CMC1.5.1}).                 											
\end{theo}

Without loss of generality we normalize
$$
H=1.
$$
The quaternion $\Phi(z,{\bar z},\Lambda)$  solving the system (\ref{R32.2.16},
\ref{R32.2.18}) with $Q^{t}=\Lambda Q$ describes the moving frame $F_z,F_{\bar z},N$ 
(\ref{R32.2.14}, \ref{R32.2.15}) of the corresponding surface. Knowing
the family $\Phi(z,{\bar z},\Lambda)$ in a neighbourhood of $\Lambda=e^{2it}$ allows us
to derive an immersion formula \cite{BoCMC} without integration the frame
with respect to $z,{\bar z}$, but just by differentiation by $t$.
Before presenting this important formula we pass to 
a gauge equivalent frame function
\begin{equation}
\Phi_0=e^{-u/4}\left(\begin{array}{cc} {1\over\sqrt{i \lambda}} &0\\
0&\sqrt{i \lambda}\end{array}\right) \Phi,\qquad \Lambda=\lambda^2,       \label{CMCT1.2}
\end{equation}
normalized by
\begin{equation}
\Phi_0(z,\bar{z},\lambda=e^{it})\in SU(2),\qquad t\in \BbbR.  \label{CMCT1.3}
\end{equation}
This lecture deals essentially with the theory of CMC tori. Since the canonical bundle in
this case is trivial, introducing a global complex coordinate, one can describe spinors
in terms of doubly-periodic functions (see Section \ref{s2.3}). 
\begin{theo}						\label{2.t2}
Let $\Phi_0(z,{\bar z},\lambda=e^{it})$ be a solution of 
the system
\begin{equation}
\Phi_{0\, z}=U_0(\lambda)\Phi_0,\qquad \Phi_{0\, \bar z}=V_0(\lambda)\Phi_0,
							\label{Phi_z}
\end{equation}
\begin{eqnarray}
U_0(\lambda)=\left(\begin{array}{cc}\displaystyle\frac{u_z}{4}  & 
i\lambda Q e^{-u/2}\\
\lambda \displaystyle \frac{i}{2} e^{u/2}  &  -\displaystyle\frac{u_z}{4}\end{array}\right),\ 
V_0(\lambda)= \left(\begin{array}{cc}-\displaystyle\frac{u_{\bar z}}{4}  & 
\displaystyle \frac{i}{2\lambda} e^{u/2}\\  
\displaystyle\frac{i}{\lambda}{\bar Q}e^{-u/2} &  
\displaystyle\frac{u_{\bar z}}{4}\end{array}\right)	\label{CMC1.5.2}
\end{eqnarray}
normalized by (\ref{CMCT1.3}).
Then $F$ and $N$, defined by the formulas
\begin{eqnarray}					\label{SymCMC}
F=-\Phi_0^{-1}\frac{\partial} {\partial t} \Phi_0 +
{i\over 2} \Phi_0 ^{-1} \sigma_3 \Phi_0,\qquad 
N=-i\Phi_0^{-1}\sigma_3\Phi_0,
\end{eqnarray}
describe a CMC surface and its Gauss map, with metric $e^u$, mean curvature 
$H=1$, and Hopf differential $Q^{t}=e^{2it} Q$.

Conversely, let $F$ be a conformal pa\-ra\-met\-ri\-za\-tion of a CMC surface
with metric $e^u$,  mean curvature $H=1$, and  Hopf differential 
$Q^{t}$. Then $F$ is given by formula (\ref{SymCMC}) where
$\Phi_0$ is a solution of (\ref{Phi_z}, \ref{CMC1.5.2}) as above.												
\end{theo}

{\it Proof}. 
First we note that both $F$ and $N$ are imaginary 
quaternions and therefore can be identified  with vectors in $\BbbR^3$. 
By identification (\ref{CMCT1.2}) the system (\ref{CMC1.5.2}) coincides with the 
quaternionic representation (\ref{R32.2.18}) for the equations for the moving frame
with the Hopf differential $\lambda Q$. Differentiating (\ref{SymCMC})
we get
\begin{eqnarray*}
F_z & = & -\Phi_0^{-1}\frac{\partial U_0(\lambda)}{\partial t}
\Phi_0 +{i\over 2}\Phi_0^{-1}[\sigma_3,U_0(\lambda)]\Phi_0
=-ie^{u/2}\Phi^{-1}\left(
\begin{array}{cc}0 & 0\\ 1 & 0 \end{array}\right)\Phi,\\
F_{\bar z} & = & -ie^{u/2}\Phi^{-1}\left(\begin{array}{cc} 0 & 1 \\ 0 & 0
\end{array}\right)\Phi,
\end{eqnarray*}
which coincides with (\ref{R32.2.15}).

{\em Remark}. In a neighborhood of a non-umbilic point $Q\not =0$ by a 
conformal change of coordinate $z\rightarrow w(z)$ one can always
normalize 
$$
Q={1\over 2}.
$$ 
Thus, umbilic free CMC surfaces are isothermic. 
In this pa\-ra\-met\-ri\-za\-tion the Gauss 
equation becomes the elliptic sinh-Gordon equation 
\begin{equation}
u_{z\bar{z}}+\sinh u=0.						\label{sinh-Gordon}
\end{equation}

\subsection{Loop group formulation}		\label{s2.2}

The matrices
$$
A=U_0+V_0,\qquad B=i(U_0-V_0)
$$
corresponding to real vector fields $\partial_x=\partial_z+\partial_{\bar{z}}$ and 
$\partial_y=i(\partial_z-\partial_{\bar{z}})$
belong to the loop algebra
$$
g_H[\lambda ] =\{ \xi:S^1\to su(2): 
\xi (-\lambda)=\sigma_3 \xi (\lambda)\sigma_3 \},
$$
and $\Phi_0$ in (\ref{Phi_z}) lies in the corresponding loop group
\begin{equation}
G_H[\lambda ] =\{ \phi:S^1\to SU(2): 
\phi (-\lambda)=\sigma_3 \phi (\lambda)\sigma_3 \}.     \label{H1.8}
\end{equation}
Here $S^1$ is the set $| \lambda |=1$.
When defined for general complex $\lambda$, elements of $g_H[\lambda ]$ and $G_H[\lambda ]$ satisfy the
real reduction
$$
\xi(\lambda)=\sigma_2\overline{\xi({1\over \bar{\lambda}})}\sigma_2,\qquad
\phi(\lambda)=\sigma_2\overline{\phi({1\over \bar{\lambda}})}\sigma_2.
$$
For applying analytic methods of the theory of integrable systems it is crucial that
CMC surfaces can be characterized in terms of this loop group completely without
referring to the special geometric nature of the coefficients of $A$ and $B$. It is not difficult 
to prove the following

\begin{theo}							\label{2.t4}
Let $\phi:D\to G_H[\lambda]$ be a smooth map on $D\subset \BbbC$ satisfying
$$
\phi_z\phi^{-1}=A\lambda +B
$$
with $A:D\to GL(2,\BbbC)$. Then the gauge equivalent
$$
\Phi_0=\exp ({i\over 2}\arg A_{21}\, \sigma_3)\phi
$$
satisfies (\ref{Phi_z}) with $U_0, V_0$ of the form (\ref{CMC1.5.2}) and describes 
the conformal  frame of the immersion
$$
F=-\phi^{-1}{\partial\over\partial t}\phi+ {i\over 2}\phi^{-1}\sigma_3 \phi,\qquad
\lambda=e^{it}
$$
of $D$ in $\BbbR^3$ with the  mean curvature $H=1$.
\end{theo}

\subsection{CMC Tori. Analytic formulation}			\label{s2.3}

Methods of Section \ref{s2.2} can be used not only in local but also in global 
studies of CMC surfaces. It is a classical result of Hopf \cite{Hop} that the only 
CMC surface of genus zero is a round sphere. Indeed the holomorphic quadratic 
differential $Q\, dz^2$ on
a sphere must vanish identically. Then (\ref{R31.2.4}) implies in particular 
$N+F=C={\rm const}$, which yields $<F-C,F-C>=1$.

Classification of CMC tori is not as simple as of spheres
but analytic tools enable us to achieve success in this case also. The reason 
for a simplification in the case $g=1$ is the fact that, unlike the case of Riemann 
surfaces of genus $g \ge 2$, on a torus it is possible to introduce a global complex 
coordinate. 
 
Any Riemann  surface of genus 1 is conformally equivalent to the factor  of 
the complex plane by a lattice $\BbbC /{\cal L}$. The corresponding conformal 
parametrization of a torus is given by a doubly-periodic mapping
$F:\BbbC /{\cal L}\to \BbbR^3$. The metric and the Hopf differential in this 
parametrization are described 
by doubly-periodic functions $u(z,\bar{z}),\  Q(z,\bar{z}).$ Note that $H={\rm const}$
implies $Q_{\bar{z}}=0$ and $Q(z)$ is a bounded elliptic function, thus a constant.
This constant is not zero, otherwise, as follows from the consideration above
the surface is a sphere. Thus CMC tori have no umbilic points. As before we normalize 
the Gauss equation to (\ref{sinh-Gordon}) by $Q={1\over 2}$.

Denoting the  generators of ${\cal L}$ by
$$Z_1=X_1+iY_1,\qquad Z_2=X_2+iY_2$$
one obtains the following
\begin{prop}
Any torus with mean curvature $H=1$ can be conformally parametrized by a 
doubly-periodic immersion $F:\BbbC\to \BbbR^3$ 
$$
F(z+Z_i,\bar{z}+\bar{Z}_i)=F(z,\bar{z}),\qquad i=1,2
$$
with the Hopf differential 
$Q={1\over 2}$. In this parametrization the metric $u(z,\bar{z})$ is a doubly-periodic 
solution to the elliptic sinh-Gordon equation (\ref{sinh-Gordon}). 
\end{prop}

Note that due to the ellipticity of equation (\ref{sinh-Gordon}) all CMC tori are real 
analytic. 

To describe all CMC tori one should solve the following problems.
\begin{enumerate}

\item[1.] Describe all doubly-periodic solutions $u(z,\bar{z})$
of the elliptic sinh-Gordon equation (\ref{sinh-Gordon}).

\item[2.] Integrate linear system (\ref{Phi_z}) with 
\begin{eqnarray}					\label{U_0,V_0}					
U_0(\lambda)={1\over 2}\left(\begin{array}{cc}\displaystyle\frac{u_z}{2}  & 
i\lambda  e^{-u/2}\\
i\lambda\displaystyle  e^{u/2}  &  -\displaystyle\frac{u_z}{2}\end{array}\right),\ 
V_0(\lambda)= {1\over 2}\left(\begin{array}{cc}-\displaystyle\frac{u_{\bar z}}{2}  & 
\displaystyle \frac{i}{\lambda} e^{u/2}\\  
\displaystyle\frac{i}{\lambda}e^{-u/2} &  
\displaystyle\frac{u_{\bar z}}{2}\end{array}\right)	\label{Lax}
\end{eqnarray}
to find $\Phi _0(z,\bar{z},\lambda =e^{it} )$.

\item[3.] Formula (\ref{SymCMC}) for $F$ describes the corresponding CMC 
immersion.
In general, this immersion is not doubly-periodic. One should specify 
parameters of the solution $u(z,\bar{z})$, which yield doubly-periodic 
$F(z,\bar{z}).$
\end{enumerate}

These three problems can be solved simultaneously using methods 
of the finite-gap integration theory. In the rest of the lecture we give an idea of 
how this solution is found. 

\subsection{Higher flows and the fundamental theorem}		\label{s2.4}

Let $u(z, \bar{z})$ be a solution of the sinh-Gordon equation. The 
perturbation
$
u_\epsilon (z, \bar{z})=u(z, \bar{z}) +\epsilon v(z, \bar{z}) 
$
of $u(z, \bar{z})$ satisfies (\ref{sinh-Gordon}) up to the terms of
order $O(\epsilon^2)$
if and only if $v(z, \bar{z})$ is a solution of the linearized elliptic
sinh-Gordon equation
\begin{equation}
(\partial_{z\bar{z}} +\cosh u(z, \bar{z})) v(z, \bar{z}) =0.  \label{CMCT2.33}
\end{equation}

The elliptic sinh-Gordon equation is one of the possible real versions of the
sine-Gordon equation, which is one of the basic models of the theory of
integrable systems. Integrable systems possess infinitely many conservation laws,
which induce infinitely many commuting flows of the corresponding dynamical system.
In particular, applying standard algebraic tools of the theory to the sine-Gordon
equation one can prove that there exists
$v(u_z,\ldots, u_z^{(k)})$,
which solves (\ref{CMCT2.33}) and is a polynomial of all its arguments. 
Such a polynomial can be treated 
as a tangential vector field to the space of solutions of the elliptic
sinh-Gordon equation. These vector fields induce flows on the phase space of
the dynamical system (\ref{sinh-Gordon}), which in the theory of solitons are called
{\em higher flows}.

There exists a regular algebraic description of these commuting flows through 
{\em formal Killing field} (see \cite{FPPS}), which is in our case a symmetric 
$K_0(-\lambda)=\sigma_3 K_0(\lambda)\sigma_3$ 
formal power series solution
\begin{equation}
K_0(\lambda)=\sum_{m=1}^\infty K_m\lambda ^{-m}          \label{CMCT2.11}
\end{equation} 
of
$$
K_0(\lambda)_z=\big[ U_0(\lambda), K_0(\lambda)\big], \qquad
K_0(\lambda)_{\bar {z}}=\big[ V_0(\lambda), K_0(\lambda)\big].     
$$
Coefficients $K_m$ can be computed recursively. 

\begin{lemm}
The diagonal terms
$$
K_{2n}= v_n\sigma_3, \qquad n=1,\ldots        
$$
of the formal Killing field (\ref{CMCT2.11}) define tangential vector fields $v_n$
\begin{equation}
(\partial_{z\bar{z}} +\cosh u) v_n=0.                \label{CMCT2.35}
\end{equation}
$v_n$ are polynomials of $u_z,\ldots, u_z^{(2n-1)}$.
\end{lemm}
Any complex vector field $v_n$ generates two real tangential vector fields
$$
w_{2n-1}=v_n+\bar{v}_n,\qquad w_{2n}=i(v_n-\bar{v}_n), \quad n=1,\ldots
$$
\begin{lemm}
Let $u(z,\bar{z}) $ be a doubly periodic solution
$$ 
u(z+Z_i,\bar{z}+\bar{Z}_i)=u(z,\bar{z}),\quad i=1,2 \quad {\rm Im}\ Z_1/Z_2\not= 0
$$
of the elliptic sinh-Gordon equation and $w_n,n=1,\ldots$ be the  corresponding
tangential real vector fields. Only finitely many
tangential vectors $w_n$ are linearly independent.
\end{lemm}
{\em Proof}.
All $w_n$ are also doubly-periodic. Equation (\ref{CMCT2.35}) determines an elliptic
operator $L$  on the torus $T:$
$$
Lw_n=(\partial_z \partial_{\bar{z}}+\cosh u)w_n=0.
$$
It is well known that the spectra of this operator is discrete, 
which implies in particular that all the eigenspaces are finite dimensional.
All tangential vectors $w_n$ belong to the kernel of $L$. 
This observation proves the lemma
$${\rm dim\ span} \{w_n\}_{n=1, \ldots}<\infty.$$

This lemma is the reason for the existence of a {\em polynomial Killing field}.
 
\begin{theo}          \label{CMCT2.t1}
Let $u(z,\bar{z})$ be a doubly-periodic solution of the elliptic sinh-Gordon equation 
(\ref{sinh-Gordon}), and $U_0, V_0$ are given by (\ref{Lax}). 
Then in the loop algebra $g_H[\lambda]$ there exists a polynomial Killing field
\begin{equation}
W_0(\lambda) =\sum _{n=-(2N-1)}^{2N-1} W_n\lambda^n       \label{CMCT2.5}
\end{equation}
which satisfies
\begin{eqnarray}
W_0 (\lambda)_z &=& \big[ U_0 (\lambda),W_0 (\lambda) \big],      \nonumber\\
W_0 (\lambda)_{\bar{z}} &=& \big[ V_0 (\lambda),W_0 (\lambda) \big].  \label{CMCT2.9}
\end{eqnarray}
\end{theo}

This fundamental theorem in different forms appeared first in \cite{PS, Hit}.

The coefficients of $W_n,n>0$ are polynomials of $u_z, \ldots, u_z^{(2N-1-n)}$
and $e^{\pm u/2}$, the coefficients of $W_n,n<0$ are polynomials of
$u_{\bar z}, \ldots, u_{\bar z}^{(2N-1+n)}$ and $e^{\pm u/2}$, $W_0$ is a polynominal 
of $u_z,u_{\bar z}, \ldots, u_z^{(2N-1)},u_{\bar z}^{(2N-1)}$.
The leading coefficient is of the form
$$
W_{2N-1}=\alpha \left(\begin{array}{cc} 0 & e^{-u/2}\\ e^{u/2} & 0
\end{array} \right),\qquad 
 0\not=\alpha\in \BbbC.         
$$

Solutions possessing polynomial Killing fields are called {\em solutions of finite type} or 
{\em finite-gap solutions}. The theory of finite-gap solutions is a well established branch 
\cite{DKN, BBEIM} of the theory of integrable equations.
Due to Theorem \ref{CMCT2.t1}, all doubly-periodic solutions of the elliptic sinh-Gordon equation
are finite-gap. 

\subsection{The spectral curve and Baker-Akhiezer function}		\label{s2.5}

Let $u(z,\bar{z})$ be a solution of the elliptic sinh-Gordon equation with the
polynomial Killing field $W_0(\lambda )$. The curve 
\begin{equation}
{\rm det} (W_0(\lambda ) -\mu I)=0                              	\label{CMCT3.2}
\end{equation}
is called the {\em spectral curve} of the solution $u(z,\bar{z})$. The spectral curve is
independent of $z, \bar{z}$.

Compactified at $\mu=\infty$ the hyperelliptic curve (\ref{CMCT3.2}) determines a 
compact Riemann surface $\hat{C}$ of genus $\hat{g}$. Due to symmetries of the loop 
algebra $g_H[\lambda ]$, besides the hyperelliptic involution $(\mu,\lambda)\to (-\mu,\lambda)$ 
it possesses two more involutions: a holomorphic
\begin{equation}
\pi: (\mu, \lambda )\to (\mu, -\lambda )                         	\label{CMCT3.10}
\end{equation}
and an anti-holomorphic $\hat{\tau}: (\mu, \lambda )\to (\bar{\mu}, {1 \over\bar{\lambda}})$.

The factor Riemann surface $C=\hat{C}/\pi$ plays central role for the explicit 
construction of section \ref{s2.6}. The covering $\hat{C}\to C,\, (\mu, \lambda )\mapsto
(\mu,\Lambda), \, \Lambda:=\lambda^2 $ is unramified and $C$ is a Riemann surface of genus $g$, 
where $\hat{g}=2g-1$. The anti-holomorphic involution 
\begin{equation}
\tau: (\mu, \Lambda )\to (\bar{\mu}, {1 \over\bar{\Lambda}}). 	\label{CMCT3.11}
\end{equation}
acts on $C$.

Due to (\ref{CMCT2.9}) the system
\begin{eqnarray*}
\phi_z=U_0\phi, \qquad \phi_{\bar{z}}=V_0\phi, \qquad W_0\phi=\mu \phi 
\end{eqnarray*}
has a common vector valued solution $\phi(P,z,\bar{z})$, which is called the 
{\em Baker-Akhiezer function}. Here $P=(\mu,\lambda)$ is a point on $\hat{C}$. 
In the finite-gap integration theory of the sine-Gordon equation, usually a gauge
equivalent function
$$
\psi= \left(\begin{array}{cc} e^{u/4} & 0\\ 0 & e^{-u/4} \end{array} \right) \phi
$$
is used. Immersion formula (\ref{SymCMC}) is obviously invariant with respect to this
transformation.

Suitably normalized 
$$
\psi(P,0,0)=\left(\begin{matrix} 1 \cr *
\end{matrix}\right)
$$ 
this function satisfies
\begin{equation}					\label{pi-symmetry}
\psi(\pi P,z,\bar{z})=\sigma_3 \psi(P,z,\bar{z}).
\end{equation}
The Baker-Akhiezer function $\psi$ has essential singularities at the points 
$\infty^\pm, 0^\pm\in\hat{C}$ defined by $\lambda(\infty^\pm)=\infty, \lambda(0^\pm)=0$. 
The involution $\pi$ interchanges these points $\pi(\infty^+)=\infty^-, \pi(0^+)=0^-$.
Denote their projection on $C$ by $\infty$ and $0$ respectively.
Due to the symmetry (\ref{pi-symmetry}) the pole divisor of $\psi$
on $\hat{C}\setminus \{\infty^\pm,0^\pm\}$ is the lift of a divisor $\cal D$ on 
$C\setminus \{\infty,0\}$.

Finally after some computations one can prove the following analytic properties of $\psi$.
\begin{theo}						\label{BAproperties}
The Baker-Akhiezer function $\psi$ possesses the following analytic properties:
\begin{enumerate}
\item[1.]
$\psi$ is transformed by (\ref{pi-symmetry}) under the action of the involution $\pi$,
\item[2.]
$\psi$ is  meromorphic on ${\hat C}\setminus
\{\infty^\pm, 0^\pm \}$. The pole divisor $\cal D$ of 
$\psi$ on $C\setminus\{\infty, 0\}$ is independent of $z, \bar{z}$, and is a non-special 
divisor of degree $g$. The Abel map ${\cal A}({\cal D})$ of $\cal D$ 
on $C$ satisfies
\begin{equation}					\label{realD}
{\cal A}({\cal D}-\tau{\cal D})={\cal A}(0-\infty),
\end{equation} 
\item[3.]
$\psi$ has essential singularities at the points $\infty^\pm, 0^\pm$ of the form
\begin{eqnarray*}
\psi(P,z,\bar{z})&=&\left(\left(\begin{matrix}1 \cr \pm 1\end{matrix}\right)+o(1)\right)
\exp(\pm{i \lambda z\over 2}),\qquad P\to\infty^\pm.\\                
\psi(P,z,\bar{z})&=& O(1)\exp(\mp{i \bar{z}\over 2 \lambda}),\qquad P\to 0^\pm.
\end{eqnarray*} 
\end{enumerate}
\end{theo}

\subsection{Baker-Akhiezer function. Formulas}			\label{s2.6}

Due to the symmetry (\ref{pi-symmetry}) the Baker-Akhiezer function $\psi$ can be
described in terms of the data $\{ C, {\cal D}\}$. Here $C$ is a hyperelliptic
Riemann surface of genus $g$ with the anti-holomorphic involution (\ref{CMCT3.11}) and 
branch points at $\lambda=0,\infty$ and $\cal D$ is a non-special divisor of 
degree $g$ on $C$ satisfying (\ref{realD}). We call these data admissible.
It is crucial that the construction of Section \ref{s2.5} can be reversed and a result
similar to Theorem \ref{2.t4} holds. 

\begin{theo}
Let $\{ C, {\cal D}\}$ be admissible data. There exists a Baker-Akhiezer function 
$\psi$ with these data and $\psi$ is uniquely characterized by the analytic 
properties listed in Theorem \ref{BAproperties}. 
\end{theo}

Admissible $\{ C, {\cal D}\}$ generate a finite-gap solution of the elliptic
sinh-Gordon equation and thus a surface with constant mean curvature $H=1$, which we call
a CMC surface of {\em finite type}. It follows from Sections \ref{s2.4}, \ref{s2.5} that 
all CMC tori are CMC surfaces of finite type.

The Baker-Akhiezer functions and hence CMC surfaces of finite type can be described 
explicitly.  Let 
$$
M^2=\Lambda\prod_{i=1}^g(\Lambda-\Lambda_i)(\Lambda-{1\over \bar{\Lambda}_i}),\qquad
\mid\Lambda_i\mid<1 \quad \forall i
$$
be a non-singular hyperelliptic curve $C$ with an anti-holomorphic involution 
$\tau:\Lambda\to{1\over\bar{\Lambda}}$. Choose a canonical homology basis 
$a_1, b_1,\ldots,a_g, b_g$ with $a_i$-cycles surrounding the cuts
$[\Lambda_i,{1\over \bar{\Lambda}_i}]$, i.e. $\tau a_i= -a_i$. Let $\omega_1,\ldots,\omega_g$  
be the dual basis 
$$
\int_{a_n}\omega_m=2\pi i\delta_{nm}
$$
of holomorphic differentials. The period matrix
$$
B_{nm}=\int_{b_n}\omega_m
$$ 
determines the Riemann theta function
$$
\theta (u)=\sum_{k\in\BbbZ^g}\exp ({1\over 2}(Bk,k)+(u,k)),\qquad u\in\BbbC^g.
$$
We need also the Abelian differentials of the second kind $\Omega_\infty, \Omega_0$ 
normalized by the condition
$$
\int_{a_n}\Omega_\infty=\int_{a_n}\Omega_0=0,\qquad  n=1,\ldots,g
$$
and the following asymptotics at the singularities
\begin{eqnarray*}
&\Omega_\infty\to d\sqrt{\Lambda},\qquad \Lambda\to\infty\\
&\Omega_0\to d{1\over\sqrt{\Lambda}},\qquad \Lambda\to 0.
\end{eqnarray*}
Denote the vector of $b$-periods of $\Omega_\infty$ by
$$
U=(U_1,\ldots,U_g),\qquad U_n=\int_{b_n}\Omega_\infty.
$$
Finally note that in explicit description one can replace the divisor $\cal D$ of admissible 
data by its Abel map $D\in Jac(C)$. One can show that in the chosen normalizations the reality 
condition (\ref{realD}) is equivalent to the condition $D\in i\BbbR^g$.

\begin{theo}							\label{BAformula}
The Baker-Akhiezer function with the data $\{ C, D\}$ is given by the formulas
\begin{eqnarray*}
\psi_1(P,z,\bar{z})&=&{\theta(\int_\infty^P\omega+W)\theta(D)\over
	\theta(\int_\infty^P\omega+D)\theta(W)}
	\exp ({i\over 2}\int_\infty^P(z\Omega_\infty+\bar{z}\Omega_0))\\
\psi_2(P,z,\bar{z})&=&{\theta(\int_\infty^P\omega+W+\Delta)\theta(D)\over
	\theta(\int_\infty^P\omega+D)\theta(W+\Delta)}
	\exp ({i\over 2}\int_\infty^P(z\Omega_\infty+\bar{z}\Omega_0)).
\end{eqnarray*}
Here $\Delta=\pi i(1,\ldots,1)$, the vector $D\in i\BbbR^g$ is arbitrary, 
$$
W=i{\rm Re}\,(Uz)+D
$$
and the integration paths in all the integrals are identical. 
The corresponding solution to the sinh-Gordon equation is given by
\begin{equation}						\label{SGformula}
u(z, \bar{z})=2\log {\theta(W+\Delta)\over \theta(W)}.
\end{equation}
\end{theo}

Applying now Theorem \ref{2.t2} to $\psi(P,z,\bar{z})$ with 
$P=P_0=(M_0,\Lambda_0),\, \mid\Lambda_0\mid=1$ one arrives at the
following final formulas for CMC immersions of finite type \cite{BoCMC}.

\begin{theo}							\label{CMCfinite-gap}
The quaternion valued solution $\Phi(z, \bar{z}, \lambda_0)$ of the linear system 
(\ref{Phi_z}, \ref{U_0,V_0}) with the finite-gap coefficient (\ref{SGformula}) is given by
$$
\Phi={i\over\sqrt{\theta(W)\theta(W+\Delta)}}
\left(\begin{array}{cc}\theta(W+l)  & \theta(W-l)\\
\theta(W+\Delta+l)  &  -\theta(W+\Delta-l)\end{array}\right)
\exp (i\sigma_3{\rm Re}\,(z L)),
$$
where $l=\int_\infty^{P_0}\omega$ is the Abel map of $P_0=(M_0, \Lambda_0)$ chosen on 
the unit circle 
$\Lambda_0=\lambda^2=e^{2it}$ and $L=\int_\infty^{P_0}\Omega_\infty$. The matrix $\Phi$ is normalized by
$$
\det \Phi= 2{\theta(l)\theta(l+\Delta)\over \theta(0)\theta(\Delta)}.
$$
The corresponding CMC immersion is given by (\ref{SymCMC}). This immersion is doubly-periodic
granted a lattice ${\cal L}$ with the basic vectors $Z_1, Z_2$ exists such that
\begin{equation}						\label{period-frame}
{\rm Re}\, (Z_k U)\in 2\pi \BbbZ^g, \quad 
{\rm Re}\, (2 Z_k \int_\infty^{P_0}\Omega_\infty) \in 2\pi \BbbZ, 
\qquad k=1,2,
\end{equation}
and the differential $\Omega_\infty$ vanishes at the point $P_0$
\begin{equation}						\label{period-immersion}
\Omega_\infty(\lambda_0)=0.
\end{equation}
\end{theo}

CMC tori are singled out from general quasiperiodic immersions of finite type by the 
periodicity conditions (\ref{period-frame}, \ref{period-immersion}), which are in fact  
conditions on the corresponding hyperelliptic curve $C$ of genus $g$ only. 
One can show \cite{Jag, BoCMC} that there are no CMC tori with $g=1$ and that for $g>1$ there 
exists a discrete set of spectral curves $C$ generating CMC tori.
The parameter 
$D\in i\BbbR^g$ remains arbitrary. So the CMC tori with $g>2$ (changes of $D$ in the plane
${\rm span}\{ {\rm Re}\, U,{\rm Im}\, U \}$ are equivalent to reparametrization of the torus) possess 
commuting deformation flows. These deformations are area preserving \cite{BoCMC}.

\subsection{Examples of CMC Tori}				\label{s2.7}


All finite-gap solutions of the sinh-Gordon equation of genus $g=1$ and $g=2$ are 
doubly-periodic. There are no CMC tori with $g=1$.
The simplest CMC tori were found by Wente \cite{Wen} and analytically studied by Abresch \cite{Abr} 
and Walter \cite{Wal}. These tori presented in Figures 1, 2 possess a 
family of plane 
curvature lines. This implies the additional symmetry $\Lambda\to 1/ \Lambda$ of the 
corresponding spectral curve of genus $g=2$. The Wente torus in Figure 1 is comprised
of three congruent fundamental domains shown in Figure 2.    

Spectral curves of genus $g=2$ without additional symmetries also generate CMC tori. An example
is presented in Figure 3.

Taking spectral curves with $g=3$ and the symmetry $\Lambda\to 1/\Lambda$ one obtains all
CMC tori with spherical curvature lines. The fundamental domain of such an example is shown in 
Figure 4. 

Figure 5 visualizes a CMC torus corresponding to a curve of genus $g=5$. This
torus possesses a 3-parameter family of area preserving deformations.
Finally Figure 6 presents classical surfaces of Delaunay which correspond to spectral
curves of genus $g=1$ and are CMC surfaces of revolution.  

  
The figures of this section are produced by Matthias Heil  using formulas presented in 
this lecture and the software for calculations on hyperelliptic Riemann surfaces developed 
by him for Sfb288 in Berlin. Further examples can be found in \cite{Hei}.

\newpage 
\section{Lecture. Bonnet Pairs}					\label{s3}

In this lecture we present some preliminary results on local and global geometry of Bonnet pairs.
This is a work in progress.
 
\subsection{Basic Facts about Bonnet Pairs}			\label{s3.1}

Let ${\cal F}_1, {\cal F}_2\subset\BbbR^3$ be a smooth Bonnet pair (Bonnet mates), i.e. two
isometric non-congruent surfaces with coinciding mean curvatures at the corresponding
points. As conformal immersions of the same Riemann surface
$$
F_1:{\cal R}\to \BbbR^3,\qquad F_2:{\cal R}\to \BbbR^3
$$
they are described by the corresponding Hopf differentials $Q_1, Q_2$, the common metric
$e^u\, dzd\bar{z}$ and the mean curvature function $H$. Since the surfaces are non-congruent
the Hopf differentials differ $Q_1\not\equiv Q_2$.

The Gauss-Codazzi equations immediately imply
\begin{prop}							\label{3.p1}
Let $Q_1$ and $Q_2$ be the Hopf differentials of a Bonnet pair $F_{1,2}\to\BbbR^3$. Then
\begin{equation}						\label{hol}
h=Q_2-Q_1
\end{equation}
is a holomorphic quadratic differential  $h\, dz^2$  on $\cal R$ and
\begin{equation}						\label{mod}
\mid Q_1 \mid=\mid Q_2 \mid.
\end{equation}
\end{prop}

Due to the second statement of Proposition \ref{3.p1} the umbilic points of ${\cal F}_1$ and
${\cal F}_2$ correspond. Denote by
$$
{\cal U}=\{ P\in{\cal R}:\ Q_k(P)=0\}
$$
the corresponding set of umbilic points on $\cal R$.

\begin{prop}							\label{3.p2}
Let $Q_1$ and $Q_2$ be the Hopf differentials of a Bonnet pair $F_{1,2}\to\BbbR^3$. Then there
exist a holomorhic quadratic differential $h$ on $\cal R$ and a smooth real valued function
$\alpha:{\cal R}\to\BbbR$ such that
\begin{equation}						\label{Q12}
Q_1={1\over 2}h (i\alpha -1),\qquad Q_2={1\over 2}h (i\alpha +1).
\end{equation}
\end{prop}

{\it Proof}. Define a smooth quadratic differential $g\, dz^2$ by
$$
g=Q_1+Q_2.
$$
Identity (\ref{mod}) implies
$$
h\bar{g}+ g\bar{h}=0.
$$
Thus the quotient
$$
\alpha=-i{g\over h}
$$
is a real valued smooth function $\alpha:{\cal R}\setminus {\cal U}_h\to\BbbR$, where
$$
{\cal U}_h=\{ P\in{\cal R}:\ h(P)=0\}
$$
is the zero set of $h$. Let us show that $\alpha$ can be smoothly extended to the whole of
$\cal R$. At any $z_0\in {\cal U}_h$ the holomorphic differential $h$ is of the form
$$
h(z)=(z-z_0)^J h_0(z),\quad h_0(z_0)\neq 0,\quad J\in\BbbN.
$$ 
 Real-valuedness of $\alpha$ near $z_0$ implies
$$
g(z)=(z-z_0)^J g_0(z)
$$
with $g_0$ smooth, which in its turn implies smoothness of $\alpha$ at $z_0$.

\begin{coro}
Umbilic points of a Bonnet pair are isolated. The umbilic set coincides with the zero set of $h$
$$
{\cal U}={\cal U}_h.
$$
\end{coro}

The number $-J$ where $J$ is defined above is called the {\em index} of the umbilic point. We
call the zero divisor
$$
D_u=(h)
$$
of $h$ the {\em umbilic divisor} of a Bonnet pair.

For compact Riemann surfaces ${\cal R}$, Propositions \ref{3.p1}, \ref{3.p2} imply the following
\begin{theo} \hspace{1cm}							\label{global}

\noindent (i) There exist no Bonnet pairs of genus $g=0$.

\noindent (ii) Bonnet pairs of genus $g=1$ have no umbilic points.

\noindent (iii) The umbilic divisor $D_u$ of a Bonnet pair of genus $g\ge 1$ is of degree 
$4g-4$ and
its class is $D_u \equiv 2K$, where $K$ is the canonical divisor.
\end{theo}

{\it Proof}. A holomorphic quadratic differential on a sphere vanishes identically $h\equiv
0$, which means $Q_1=Q_2$, and the surfaces are congruent. A holomorphic quadratic differential 
on a torus does not have zeros, thus ${\cal U}=\emptyset$ for tori.
\vspace{0.5cm}

The point (i) of Theorem \ref{global} was proven in \cite{LawTri}.

Taking into account the similarity of the analytic description of Bonnet surfaces and CMC surfaces
and the progress in the investigation of CMC surfaces achieved by methods from the theory of
integrable systems (see Section \ref{s2}), the most promising open problem to attack by 
these methods seems to be the problem of existence and description of Bonnet tori mates.

For tori one has ${\cal R}=\BbbC/{\cal L}$. Scaling the lattice ${\cal L}$ appropriately one can always
normalize $h=-i$, i.e.
\begin{equation}						\label{barQ}
Q_1={1\over 2}(\alpha +i),\qquad Q_2={1\over 2}(\alpha -i).
\end{equation}
The corresponding Gauss-Codazzi equations of Bonnet mates become
\begin{eqnarray}
2u_{z\bar{z}}+H^2 e^u-(1+\alpha^2)e^{-u}&=&0,		\nonumber\\
\alpha_{\bar{z}}-e^u H_z&=&0,				\label{PairsGC}\\
\alpha_{z}-e^u H_{\bar{z}}&=&0.				\nonumber
\end{eqnarray}
Note that the Gauss-Codazzi equations of isothermic surfaces $Q=\alpha/2\in\BbbR$ differ 
only slightly 
\begin{eqnarray}
2u_{z\bar{z}}+H^2 e^u-\alpha^2 e^{-u}&=&0,		\nonumber\\
\alpha_{\bar{z}}-e^u H_z&=&0,				\label{IsothGC}\\
\alpha_{z}-e^u H_{\bar{z}}&=&0.				\nonumber
\end{eqnarray}

\subsection{Lax Representation and Connection to Isothermic Surfaces}\label{s3.2}

Looking for a Lax representation for Bonnet pairs, it is natural to try to merge the frame 
equations of two Bonnet mates. For tori, cylinders, or simply connected domains it is
enough to consider the case of ${\cal R}$ being a domain $D$ in $\BbbC$. Cylinders and tori are
distinguished by the corresponding periodicity lattices $\cal L$. 
Since our main interest lies in the investigation of tori let us restrict ourselves to the case of 
umbilic free Bonnet pairs. 
As in Section \ref{s2}, introducing a global complex variable $z$ on $D$ we normalize the 
corresponding frame matrices traceless and the Hopf differentials as in (\ref{barQ}). 
The following theorem can be checked directly.

\begin{theo}							\label{PairsLax}
Normalized by (\ref{barQ}), conformal frames $\Phi_1, \Phi_2 : D\to SU(2)$ of a Bonnet pair  
\begin{eqnarray*}
&\Phi_{k\, z}=U_k \Phi_k,\qquad \Phi_{k\, \bar{z}}=V_k \Phi_k\\
&U_k=\left(\begin{array}{cc}{u_z\over 4} & - Q_k e^{-u/2}\\
{H\over 2}e^{u/2} & -{u_z\over 4}\end{array}\right),\qquad
V_k=\left(\begin{array}{cc} -{u_{\bar{z}}\over 4} & -{H\over 2}e^{u/2}\\
 \bar{Q}_k e^{-u/2} & {u_{\bar{z}}\over 4}\end{array}\right)
\end{eqnarray*}
can be extended 
\begin{equation}						\label{Phi(0)}
\Phi(z, \bar{z}, \lambda=0)=
\left(\begin{array}{cc} \Phi_1 & 0\\
0 & \Phi_2 \end{array}\right)(z, \bar{z}) 
\end{equation}
to $\Phi(z,\bar{z},\lambda)$ satisfying 
\begin{equation}						\label{PairsPhi}
\Phi_{z}=U \Phi,\qquad \Phi_{\bar{z}}=V \Phi
\end{equation}
with 
\begin{eqnarray}
U&=&\left(\begin{array}{cc} U_1 & 
\begin{array}{cc} 0 & 0 \\
-i\lambda e^{u/2} & 0\end{array}\\
\begin{array}{cc} 0 & -i\lambda e^{-u/2} \\
0 & 0\end{array} & U_2
\end{array}\right), 						\label{PairsUV}\\								
V&=&\left(\begin{array}{cc} V_1 & 
\begin{array}{cc} 0 & -i\lambda e^{u/2} \\
0 & 0\end{array}\\
\begin{array}{cc} 0 & 0 \\
-i\lambda e^{-u/2} & 0\end{array} & V_2
\end{array}\right).						\nonumber
\end{eqnarray}
Conversely, the linear system (\ref{PairsUV}) with $Q_1=Q\not=Q_2=\bar{Q}$
is compatible if and only if the metric $e^u$, the mean curvature function $H$, and the Hopf 
differentials $Q_1, Q_2$ on $D$ satisfy the Gauss-Codazzi equations of Bonnet pairs.
Conformal frames of the Bonnet mates are determined through (\ref{Phi(0)}) by a suitably 
normalized common solution $\Phi(z,\bar{z},\lambda)$ of (\ref{PairsPhi}) evaluated at 
$\lambda=0$.   
\end{theo} 

{\it Remark}. 
In the general case of an arbitrary Riemann surface $\cal R$ and holomorphic quadratic 
differential $h$, a spinor form of the Lax representation for Bonnet pairs similar to 
(\ref{PairsLax}) can be easily made by merging the spinor frames (see Section \ref{s1}) of the
corresponding surfaces.
\vspace{0.5cm}
 
{\it Remark}. In the case $Q_1=Q_2=Q=\bar{Q}$ system (\ref{PairsPhi}, \ref{PairsUV}) 
becomes a  Lax representation for isothermic surfaces in $\BbbR^3$ in isothermic 
coordinates.
\vspace{0.5cm}

The matrices $U+V$ and $i(U-V)$ corresponding to real vector fields $\partial_x$ and $\partial_y$
possess the symmetries
\begin{eqnarray}
A(-\lambda)&=& \left(\begin{array}{cc} {\bf 1} & 0\\
0 & -{\bf 1} \end{array}\right)
A(\lambda)\left(\begin{array}{cc} {\bf 1} & 0\\
0 & -{\bf 1} \end{array}\right),					\label{-lambda}\\
\overline{A(\bar{\lambda})}&=& \left(\begin{array}{cc} \sigma_2 & 0\\
0 & \sigma_2 \end{array}\right)
A(\lambda)\left(\begin{array}{cc} \sigma_2 & 0\\
0 & \sigma_2 \end{array}\right).					\label{bar_lambda}
\end{eqnarray}
Denote by
$$
g_B[\lambda]=\{A:\BbbR\to gl(2, \BbbH)
\mid A(-\lambda)= \left(\begin{array}{cc} {\bf 1} & 0\\
0 & -{\bf 1} \end{array}\right)
A(\lambda)\left(\begin{array}{cc} {\bf 1} & 0\\
0 & -{\bf 1} \end{array}\right) \}
$$
the corresponding loop algebra, and by 
$$
G_B[\lambda]=\{\phi:\BbbR\to GL(2, \BbbH)
\mid \phi(-\lambda)= \left(\begin{array}{cc} {\bf 1} & 0\\
0 & -{\bf 1} \end{array}\right)
\phi(\lambda)\left(\begin{array}{cc} {\bf 1} & 0\\
0 & -{\bf 1} \end{array}\right) \}
$$
the corresponding loop group.

By the normalization (\ref{Phi(0)}) the solution $\Phi(z,\bar{z},\lambda)$ is
determined uniquely up to right multiplication by a matrix depending only on $\lambda$ 
\begin{equation}						\label{G(lambda)}
\Phi(z,\bar{z},\lambda)\to\Phi(z,\bar{z},\lambda)G(\lambda),\qquad G(0)=1. 
\end{equation}
$\Phi(z,\bar{z},\lambda)$ can be chosen to lie in $G_B[\lambda]$. Then the matrix
\begin{equation}						\label{Sym_Isoth}
\left(\begin{array}{cc} 0 & S\\
T & 0 \end{array}\right):=\Phi^{-1}\Phi_\lambda \mid_{\lambda=0}
\end{equation}
is off-diagonal and its coefficients are quaternion valued functions of $z$ and $\bar{z}$. 

For a description of the geometry of $S$ and $T$, the notion of isothermic surfaces in $\BbbR^4$ and of 
the dual isothermic surface is required. 
An immersion $f:D\to \BbbR^4$ is called {\em isothermic} if it is conformal and the vector $f_{xy}$ lies
in the tangent plane $f_{xy}\in\, {\rm span}\{ f_x, f_y\}$. It is convenient to describe isothermic
surfaces in $\BbbR^4=\BbbH$ in quaternionic form, i.e. as mappings $f:D\to \BbbH$ with the 
coordinates
$$
f=f_0 {\bf 1}+ f_1 {\bf i}+f_2 {\bf j}+f_3 {\bf k}.
$$
Its differential is $df=f_x\, dx+f_y\, dy$. 
An important property of an isothermic immersion $f:D\to \BbbH$ is the closedness of the form
\begin{equation}						\label{dual-isoth}
df^*:=-f_x^{-1} dx + f_y^{-1} dy
\end{equation}
The corresponding immersion determined up to translation by this form  and denoted by 
$f^*:D\to \BbbH$ is also isothermic and is called the {\em dual isothermic surface}. Note that the dual 
isothermic surface is defined through one-forms and therefore the periodicity properties of $f$ are 
not respected. The relation (\ref{dual-isoth}) is an involution. Moreover one can check that

\begin{lemm}						\label{l.dual->isoth}
The transformation (\ref{dual-isoth}) is characteristic for isothermic surfaces.
\end{lemm} 

\begin{prop}\hspace{1cm}				\label{p.Pairs-isoS3}

\noindent (i) Let $Q_1$ and $Q_2$, normalized by (\ref{barQ}), be Hopf differentials of a Bonnet 
pair $F_{1,2}:D\to {\rm Im}\, \BbbH$. Then $T:D\to \BbbH$ defined by (\ref{Sym_Isoth}) is an
isothermic surface in the three-dimensional sphere $S^3\subset\BbbR^4=\BbbH$
 and $S:D\to \BbbH$ in (\ref{Sym_Isoth}) is its dual $S=T^*$. The isothermic
surfaces $S$ and $T$ are related to the Bonnet pair by
\begin{eqnarray}
dF_1 &=& dS\, T=dT^*\, T						\label{Pairs-isoS3}\\
dF_2 &=& T\, dS=T\, dT^*.						\nonumber
\end{eqnarray}

\noindent (ii) Let $Q_1=Q_2=Q$ in (\ref{PairsUV}) be real. Then $S:D\to {\rm Im}\, \BbbH$ given 
by (\ref{Sym_Isoth}) is the isothermic surface determined by the fundamental forms with the
coefficients $e^u, H, Q$, and $T:D\to {\rm Im}\, \BbbH$ is its dual $T=S^*$.
\end{prop}

{\it Proof}. Let us prove the first statement. Formula (\ref{Sym_Isoth}) implies 
$$
d\left(\begin{array}{cc} 0 & S\\
T & 0 \end{array}\right)=\Phi^{-1}
\left( U_{\lambda} dz+V_{\lambda} d\bar{z}\right)\Phi \mid_{\lambda=0}
$$
or equivalently
\begin{equation}					\label{dSdT}
dS = e^{u/2}\Phi_1^{-1}({\bf i}dx + {\bf j}dy)\Phi_2,
dT = e^{-u/2}\Phi_2^{-1}({\bf i}dx - {\bf j}dy)\Phi_1.
\end{equation}
These frames are obviously conformal and are related by (\ref{dual-isoth}). Lemma
\ref{l.dual->isoth} implies that $T$ is isothermic. 
Moreover $dT$ can be integrated explicitly
\begin{equation}					\label{T}
T=\Phi_2^{-1}\Phi_1.
\end{equation}
Indeed, differentiating the last expression one obtains
\begin{eqnarray*}
& dT=\Phi_2^{-1}(d\Phi_1 \Phi_1^{-1}-d\Phi_2 \Phi_2^{-1})\Phi_1=\\
& e^{-u/2}\Phi_2^{-1}\left((Q_2-Q_1)\left(\begin{array}{cc} 0 & 1\\
0 & 0 \end{array}\right)dz +(\bar{Q}_1-\bar{Q}_2)\left(\begin{array}{cc} 0 & 0\\
1 & 0 \end{array}\right)d\bar{z}\right)\Phi_1
\end{eqnarray*}
which coincides with the previous expression for $dT$. Integrating, one obtains 
$T=\Phi_2^{-1}\Phi_1 +{\rm const}$. The constant can be normalized to zero by 
transformation (\ref{G(lambda)}) with an appropriate $G(\lambda)\in G_{B}[\lambda]$.
Obviously the surface given by (\ref{T}) lies in the three sphere $S^3=\BbbH_1$. 
Using (\ref{T}) we obtain 
\begin{eqnarray*}
dS\, T&=&-ie^{-u/2}\Phi_1^{-1}\left(\begin{array}{cc} 0 & d\bar{z}\\
dz & 0 \end{array}\right) \Phi_1	\\
T dS&=&-ie^{-u/2}\Phi_2^{-1}\left(\begin{array}{cc} 0 & d\bar{z}\\
dz & 0 \end{array}\right)\Phi_2
\end{eqnarray*}
which coincides with (\ref{Pairs-isoS3}).

The proof of the second claim is even simpler (see \cite{BP}). 

The next theorem is essentially due to Bianchi \cite{BianchiBonnet}. A modern version of it in terms of
quaternions is derived in \cite{KPP}.

\begin{theo}						\label{Bi-KPP}
$F_{1,2}:D\to {\rm Im}\, \BbbH=\BbbR^3$ build a Bonnet pair if and only if there exists an
isothermic surface $T:D\to \BbbH_1=S^3\subset \BbbR^4$ (or equivalently an isothermic surface
$R:D\to {\rm Im}\,\BbbH=\BbbR^3$) such that
\begin{eqnarray}
dF_1&=&dT^*\, T={1\over 2}(1-R)dR^*(1+R)			\nonumber\\
dF_2&=&T dT^*={1\over 2}(1+R)dR^*(1-R).				\label{FTR}
\end{eqnarray}
\end{theo}

{\it Proof}. Let us show first the equivalence of the representations in terms of $T$ and $R$. 
The class of isothermic surfaces is invariant under M\"{o}bius transformations. In
particular, isothermic surfaces in $S^3$ and in $\BbbR^3$ are related by stereographic
projection, which in quaternionic form  can be represented by
$$
T={1+R\over 1-R}
$$
with $R\in {\rm Im}\, \BbbH =\BbbR^3, T\in \BbbH_1=S^3$.
For the frames this implies
\begin{eqnarray*}
dT&=&2(1-R)^{-1}dR(1-R)^{-1}\\
dT^*&=&{1\over 2}(1-R)dR^*(1-R),
\end{eqnarray*}
which proves the equivalence of the two representations of the theorem.
The passing from Bonnet pairs to isothermic surfaces is proven in 
Proposition \ref{p.Pairs-isoS3}. Conversely, given $R$, the representation (\ref{FTR}) 
shows that $dF_{1,2}$ are conformal and lie in ${\rm Im}\, \BbbH$. Due to
$dF_2=T dF_1 T^{-1}$ the immersions are isometric.

Similarly one can show (see \cite{KPP}) that the forms $dF_1, dF_2$ defined by (\ref{FTR}) are
closed and that the mean curvature functions of the corresponding surfaces coincide.

To study a global version of Theorem \ref{Bi-KPP} in the case of Bonnet tori mates
let $F_{1,2}:\BbbC\to \BbbR^3$ be a Bonnet pair with doubly periodic frames $dF_{1,2}$ with the 
same period lattice ${\cal L}$. The frames $\Phi_{1,2}$ are periodic up to a sign
$$
\Phi_k(z+Z_i)=(-1)^{p_{ik}}\Phi_k(z),
$$
where $Z_1, Z_2$ are generators of ${\cal L}$ and $p_{ik}\in \BbbZ_2$ characterize the spin
structures of the immersions. The isothermic surface given by (\ref{T}) is a torus in $S^3$ if
the spin structures of $F_1$ and $F_2$ coincide ($p_{i1}=p_{i2}, i=1,2$) and is a torus in 
three dimensional real projective space $\BbbR P^3=S^3/\{-1\}$ if the spin structures differ. 
Conversely, an isothermic immersion $T:\BbbC/{\cal L}\to S^3/\{-1\}$ generates by (\ref{FTR}) 
a Bonnet pair with frames $dF_{1,2}$ defined on $\BbbC/{\cal L}$. The spin structures of these 
two surfaces are the same iff $T$ is an immersion to $S^3$ with the lattice ${\cal L}$.

\begin{coro}
Bonnet mates with doubly periodic frames $dF_{1,2}$ are in one to one
correspondence with isothermic tori in $\BbbR P^3$.
Bonnet mates with doubly periodic frames $dF_{1,2}$ and coinciding spin structure are in one 
to one correspondence with isothermic tori in $S^3$.
The corresponding relations are given by formulas (\ref{FTR}).
\end{coro}

Formula (\ref{FTR}) allows us to control the periodicity of the frame of a Bonnet pair. To be
able to control the periodicity of the immersion one needs an analog of formula (\ref{SymCMC})
describing the corresponding immersion without integration. We call a solution 
\begin{eqnarray*}
\Phi:D\times \BbbR&\to& G_B[\lambda]\\
	(z,\lambda)&\mapsto& \Phi(z,\bar{z},\lambda)
\end{eqnarray*}
of (\ref{PairsPhi}, \ref{PairsUV}) {\em normalized} (see proof of Proposition \ref{p.Pairs-isoS3}) 
if the coefficient $T$ of its decomposition 
(\ref{Sym_Isoth}) at $\lambda=0$ is a unitary quaternion $T:D\to \BbbH_1$, i.e. (\ref{T}) holds.
Obviously this solution can be extended to all $\lambda\in\BbbC$.
\begin{theo}							\label{Sym-Pairs}
Let $\Phi(z,\bar{z},\lambda)$ a normalized  solution of (\ref{PairsUV}) with $Q$ of 
the form (\ref{barQ}). Then the corresponding Bonnet pair $F_{1,2}(z,\bar{z})$ is restored by the
following coefficients of quaternionic 2 by 2 matrices
\begin{eqnarray}						
\left(\begin{array}{cc} F_1 & 0\\
0 & * \end{array}\right)&=&
{1\over 2}\Phi^{-1}\Phi_{\lambda\lambda}\mid_{\lambda=0},	\label{ImmersionPairs}\\
\left(\begin{array}{cc} * & 0\\
0 & F_2 \end{array}\right)&=&
{1\over 2}\Phi^{-1}\Phi_{\lambda\lambda}-
\left(\Phi^{-1}\Phi_\lambda \right)_\lambda {\mid_{\lambda=0}}.
\end{eqnarray}						\nonumber
\end{theo}

{\it Proof}. To prove the formula for $dF_1$ let us differentiate it by $z$ and $\bar{z}$
$$
d({1\over 2}\Phi^{-1}\Phi_{\lambda\lambda})=
\Phi^{-1} (U_\lambda dz + V_\lambda d\bar{z})\Phi_\lambda,
$$
where we used $U_{\lambda\lambda}=V_{\lambda\lambda}=0$.
Evaluating this expression at $\lambda=0$ using (\ref{Phi(0)}), (\ref{T}) and (\ref{dSdT}), one
finally obtains
$$
d({1\over 2}\Phi^{-1}\Phi_{\lambda\lambda})_{\mid_{\lambda=0}}=
\left(\begin{array}{cc} dS\, T & 0\\
0 & dT\, S \end{array}\right).
$$
Now the formula for $dF_1$ follows from (\ref{Pairs-isoS3}).
The formula for $dF_2$ is proven by an analogous computation
$$
d({1\over 2}\Phi^{-1}\Phi_{\lambda\lambda}-
\left(\Phi^{-1}\Phi_\lambda \right)_\lambda)_{\mid_{\lambda=0}}=
\left(\begin{array}{cc} S\, dT & 0\\
0 & T\, dS \end{array}\right).
$$

\subsection{Loop Group Description}				\label{s3.3}

As we have seen already in Proposition \ref{p.Pairs-isoS3} the theory developed in 
Section \ref{s3.2} includes two different cases: the case of Bonnet pairs when the Hopf 
differentials $Q_1=Q$ and $Q_2=\bar{Q}$ are different and thus generate two non-congruent 
surfaces, and the case $Q=\bar{Q}$ of isothermic surfaces in $\BbbR^3$.
The loop group $G_B[\lambda]$ and the loop algebra $g_B[\lambda]$ of Bonnet pairs are described
in Section \ref{s3.2}. In the case of isothermic surfaces $Q=\bar{Q}$, the corresponding 
loop group and algebra  are specialized further as follows:
\begin{eqnarray*} 
G_I[\lambda]&=&\{\phi\in G_B[\lambda]
\mid \phi^T(\lambda) \left(\begin{array}{cc} 0 & \sigma_2\\
\sigma_2 & 0 \end{array}\right)\phi(\lambda)=\left(\begin{array}{cc} 0 & \sigma_2\\
\sigma_2 & 0 \end{array}\right) \}\\
g_I[\lambda]&=&\{A\in g_B[\lambda]
\mid A^T(\lambda)= - \left(\begin{array}{cc} 0 & \sigma_2\\
\sigma_2 & 0 \end{array}\right)A(\lambda)\left(\begin{array}{cc} 0 & \sigma_2\\
\sigma_2 & 0 \end{array}\right) \}.
\end{eqnarray*}

As in the previous sections the main strategy for applying analytic methods of the theory of 
integrable systems consists of two steps: first, to characterize the frame equations through analytic
properties of $\Phi$ as a function of $\lambda$ without referring to the special geometric nature of the
coefficients of the frame equations and, second, to construct those $\Phi(\lambda)$ explicitly.
For this purpose it is more convenient to pass to a gauge equivalent 
\begin{equation}						\label{PhiPhi_0}
\Phi_0=\left(\begin{array}{cc} e^{-u/4}{\bf 1} & 0\\
0 & e^{u/4}{\bf 1} \end{array}\right)\Phi
\end{equation}
to (\ref{PairsPhi}), linear problem 
\begin{equation}						\label{PairsPhi_0}
\Phi_{0\, z}=U \Phi_0,\qquad \Phi_{0\, \bar{z}}=V \Phi_0
\end{equation} 
\begin{eqnarray}
U_0&=&\left(\begin{array}{cccc} 
0		& -Qe^{-u/2} 	& 0 		& 0 \\
{H\over 2}e^{u/2}&-{u_z \over 2}& -i\lambda 	& 0\\
 0 		& -i\lambda 	&{u_z \over 2}	&-\bar{Q}e^{-u/2}\\
0 		& 0		&{H\over 2}e^{u/2}& 0	
\end{array}\right), 						\label{PairsU_0V_0}\\								
V_0&=&\left(\begin{array}{cccc}
-{u_{\bar{z}} \over 2}&-{H\over 2}e^{u/2}	& 0 		& -i\lambda \\
\bar{Q}e^{-u/2}	& 0		& 0	 	& 0\\
 0 		& 0 		&0		&-{H\over 2}e^{u/2}\\
-i\lambda 	& 0		&Qe^{-u/2}	& {u_{\bar{z}} \over 2}	
\end{array}\right),						\nonumber
\end{eqnarray}
normalizing the solution at $\lambda=\infty$.
Note that all immersion formulas (\ref{ImmersionPairs}, \ref{Sym_Isoth}) are preserved under 
this gauge transformation.

\begin{theo}
Let
\begin{eqnarray*}
\Phi_0:D\times\BbbR&\to& G_B[\lambda]\\
z,\lambda &\mapsto & \Phi_0(z, \bar{z},\lambda)
\end{eqnarray*}
be a smooth mapping  satisfying 
$$
\Phi_{0\, z}\Phi_0^{-1}=A(z, \bar{z})\lambda + B(z, \bar{z})
$$
with some $A(z, \bar{z})$ and $B(z, \bar{z})$, 
having the asymptotics 
\begin{equation}					\label{Phi_0-asymp}
\Phi_0(z, \bar{z},\lambda)=(L+M(z, \bar{z}) \lambda^{-1}+o(\lambda^{-1}))
\exp (-i\lambda (J_1 z+ J_2\bar{z})) C(\lambda)
\end{equation}
at $\lambda\to\infty$ with
$$
L=\left(\begin{array}{cccc}
1	& 0	& 0 	& 1 \\
0	& 1	& 1	& 0 \\
0 	& 1 	&-1	& 0 \\
1	& 0	& 0	& -1	
\end{array}\right), \
J_1=\left(\begin{array}{cccc}
0	& 0	& 0 	& 0 \\
0	& 1	& 0	& 0 \\
0 	& 0 	&-1	& 0 \\
0	& 0	& 0	& 0	
\end{array}\right), \
J_2=\left(\begin{array}{cccc}
1	& 0	& 0 	& 0 \\
0	& 0	& 0	& 0 \\
0 	& 0 	& 0	& 0 \\
0	& 0	& 0	& -1	
\end{array}\right)
$$ 
and some $M(z, \bar{z})$ and $C(\lambda)$. If, for the coefficients of $M$, it holds that
\begin{equation}					\label{H-constraint}
M_{31}(z, \bar{z})=-M_{42}(z, \bar{z}),
\end{equation}
then 
$$
U_0:= \Phi_{0\, z}\Phi_0^{-1},\qquad V_0:= \Phi_{0\, \bar{z}}\Phi_0^{-1}
$$
can be parametrized as in (\ref{PairsU_0V_0}) with some real valued 
$u(z, \bar{z}), H(z, \bar{z})$.
If in addition
\begin{equation}					\label{I-constraint}
M_{21}(z, \bar{z})=-M_{12}(z, \bar{z}),
\end{equation}
then $Q(z,\bar{z})$ in (\ref{PairsU_0V_0}) is also real valued $Q=\bar{Q}$ and $\Phi$ 
related to $\Phi_0$ by (\ref{PhiPhi_0}) can be chosen in $G_I[\lambda]$.
\end{theo}

{\it Proof}. Due to the assumption of the theorem, $U_0$ is of the form $U_0=A\lambda +B$.
Substituting the asymptotics (\ref{Phi_0-asymp}) one obtains for the coefficients
$$
A=-iL J_1 L^{-1},\qquad B = [M L^{-1}, A].
$$
Similarly the symmetry (\ref{-lambda}) and the asymptotics (\ref{Phi_0-asymp}) imply  
$V_0=C\lambda +D$ with
$$
C=-iL J_2 L^{-1},\qquad D = [M L^{-1}, C].
$$ 
The matrices $A$ and $C$ are of the required form. For the matrices $B$ and $D$ a simple computation 
gives 
\begin{eqnarray*}
B&=&-i\left(\begin{array}{cccc}
0	& M_{12}& 0 	& 0 \\
-M_{31}	& M_{22}-M_{32}	& 0	& 0 \\
0 	& 0 	& M_{32}-M_{22}	& -M_{21} \\
0	& 0	& M_{42}	& 0	
\end{array}\right)\\
D&=&-i\left(\begin{array}{cccc}
M_{11}-M_{41}	& -M_{42}& 0 	& 0 \\
M_{21}	& 0	& 0	& 0 \\
0 	& 0 	& 0	& M_{31} \\
0	& 0	& -M_{12}	& M_{41}-M_{11}	
\end{array}\right),
\end{eqnarray*}
where we have used the symmetry of $M$
$$
\left(\begin{array}{cccc}
1	& 0	& 0 	& 0 \\
0	& 1	& 0	& 0 \\
0 	& 0 	& -1	& 0 \\
0	& 0	& 0	& -1	
\end{array}\right) M
\left(\begin{array}{cccc}
0	& 0	& 0 	& 1 \\
0	& 0	& 1	& 0 \\
0 	& 1 	& 0	& 0 \\
1	& 0	& 0	& 0	
\end{array}\right)=-M.
$$
The anti-holomorphic involution $\lambda\to \bar{\lambda}$ of the loop group implies
$$
M_{12}=\bar{M}_{21},\ M_{42}=\bar{M}_{31},\  M_{11}=-\bar{M}_{22},\ M_{32}=-\bar{M}_{41}.
$$
Finally, comparing the coefficients and using the compatibility conditions one shows that $U_0,
V_0$ are of the form (\ref{PairsU_0V_0}).

\subsection{Surfaces of Finite Type}		\label{s3.4}

We have shown that isothermic surfaces and Bonnet pairs can be studied in frames of the theory of
integrable systems. Applying standard methods one can define for these surfaces higher flows, 
the B\"{a}cklund-Darboux transformations\footnote{The B\"{a}cklund-Darboux transformations of 
isothermic surfaces are known already in local
differential geometry \cite{Dar, BianchiIsoth} and in theory of solitons \cite{Cie}}, 
finite-gap solutions etc.

In contrast to the case of CMC surfaces, the linearizations of the Gauss-Codazzi equations of 
isothermic surfaces (\ref{IsothGC}) and of Bonnet pairs (\ref{PairsGC}) are not elliptic. This fact 
prevents us from applying the arguments of Section \ref{s2} and claiming that all corresponding 
immersions with doubly periodic fundamental forms are generated by finite-gap solutions. 
We call surfaces corresponding to the finite-gap solutions of the Gauss-Codazzi equations or
equivalently surfaces with a polynomial Killing field (see Section \ref{s2}) surfaces
of {\em finite type}. This class of surfaces is worth studying, in particular, since it may
contain Bonnet tori mates. 

A polynomial Killing field $W(\lambda)$ of an isothermic surface or a Bonnet pair of finite type is 
an element of the loop algebra $g_B[\lambda]$. The coresponding spectral curve $\hat{C}$
$$
\det (W(\lambda)-\mu I)=0
$$
is a four-sheeted covering $\hat{C}\to \bar{\BbbC}\ni\mu$. Due to 
symmetries of the loop algebra it possesses a holomorphic involution 
$\pi:(\mu,\lambda)\to(\mu,-\lambda)$ and
an anti-holomorphic involution 
$\hat{\tau}:(\mu,\lambda)\to(\bar{\mu},\bar{\lambda})$. The Riemann
surface $C=\hat{C}/\pi$ is an algebraic curve of $\mu$ and $\Lambda=\lambda^2$. It also possesses an
anti-holomorphic involution $\tau:(\mu,\Lambda)\to(\bar{\mu},\bar{\Lambda})$. In the case of
isothermic surfaces in $\BbbR^3$ the Killing field lies in the subalgebra $g_I[\lambda]$, 
which implies the additional symmetry $\sigma:(\mu, \Lambda)\to (-\mu, \Lambda)$ of $C$.
The factor curve $C_{\sigma}:=C/\sigma$ is quadratic in $\mu^2$ and thus a 
hyperelliptic one. This fact simplifies the theory and leads to a more effective description 
of the corresponding surfaces.

\begin{prop}						\label{p.SCIsoth}
The spectral curve $C$ of an isothermic surface is a double covering of a hyperelliptic
curve.
\end{prop}

Let us indicate the steps of construction of the finite-gap solutions in this case.
Although the considerations are similar to those of Section \ref{s2}, technically they are more  
involved. By more or less standard technique one describes explicitly finite-gap solutions of 
the complexified system (\ref{IsothGC}), i.e. of the system corresponding to the loop algebra 
without the real reduction (\ref{bar_lambda}). The spectral curve in this case remains to be 
a double covering of a hyperelliptic curve, but does not necessarily possess the real involution 
$\tau$. 
The differentials analogous to the differentials $\Omega_1$ and $\Omega_2$ turn out to be 
Abelian differentials of the second kind with singularities at the branch points of the covering 
$C\to C_\sigma$. They are odd with respect to the involution $\sigma$. Their vectors of $b$-periods 
describing the velocity of the linear dynamics on the Jacobian lie in the odd part of the Jacobian with
respect to the involution $\sigma$. The dynamic of the corresponding nonlinear system lies 
on the Prym variety $Prym_\sigma (C)$ of the covering $C\to C_{\sigma}$.

The real reduction (\ref{bar_lambda}) leads to constraints on parameters of the finite-gap
solution. It is a technical but rather involved problem to classify all possible cases leading to real
finite-gap solutions of (\ref{IsothGC}) and thus to isothermic surfaces of finite type. This is not yet
done. Even more important is to realize this program for the system (\ref{PairsGC}). 
Note that due to explicit formulas (\ref{Sym_Isoth}, \ref{ImmersionPairs}) for the corresponding
immersions, the isolation of tori from general surfaces of finite type is then straightforward.

\newpage 
\section{Lecture. Bonnet Families}				\label{s4}

Here we classify Bonnet families.
The content of this lecture is based on results obtained jointly with Ulrich
Eitner \cite{BE1, BE2}.

\subsection{Definition of Bonnet Surfaces and Simplest Properties}\label{s4.1}

Let ${\cal F}$ be  a smooth surface in $\BbbR^3$  with non-constant mean
curvature function. ${\cal F}$ is called a {\em Bonnet surface} if it 
possesses a one-parameter family 
$$
{\cal F}_\tau, \qquad \tau\in(-\epsilon,\epsilon),\,\epsilon > 0,\ {\cal F}_0 = {\cal F}
$$
of non-trivial\footnote{We call an isometry of a surface non-trivial if it is
not induced by an isometry of the ambient space.} isometric deformations
preserving the mean curvature function.
The family 
\linebreak $({\cal F}_{\tau})_{\tau\in(-\epsilon,\epsilon)}$ is called a 
{\em Bonnet family}. A Bonnet family can be described as a conformal mapping 
\begin{equation}                                        \label{BonnetImmersion}
\begin{matrix} 
        F: &(-\epsilon,\epsilon)\times {\cal R} &  
        \longrightarrow & \BbbR^3\cr  
      	& (\tau,z)       &\mapsto  & F(\tau,z,\bar z)
\end{matrix}\qquad\epsilon >0,
\end{equation}
where $z$ is a local coordinate $z$ on the Riemann surface $\cal R$.

The set ${\cal U}\subset {\cal R}$ of preimages of the umbilic points 
on ${\cal F}_\tau$ is independent of $\tau$ (see Section \ref{s3}). Obviously, the set
$$
{\cal V} = \{ P\in {\cal R}\ :\ dH(P)=0\}
$$
of preimages of critical points of the mean curvature function on 
${\cal F}_\tau$ is also $\tau$-independent.

Similarly to Proposition \ref{3.p2} one can prove

\begin{prop}                             	\label{4.p1}
The holomorphic quadratic differential 
\begin{equation}			\label{varphi-dz^2}
        \varphi(z,\tau)\,d z^2 := 
        \dfrac{\partial}{\partial\tau}\,Q(z,\bar z,\tau) \, d z^2.
\end{equation}
vanishes exactly at the umbilic points. 
\end{prop}

In a neighbourhood of a non-umbilic point $Q(P)\neq 0$
there exist smooth real-valued functions 
$\psi(z,\bar z,\tau)$, $q(z,\bar z)$ such that 
$Q(z,\bar z,\tau) = e^{{\rm i}\,\psi(z,\bar z,\tau)}\,q(z,\bar z)$.
Differentiating we obtain 
\begin{equation}                        \label{isothermicQ}
        Q(z,\bar z,\tau) = -{\rm i}\,
        \dfrac{\varphi(z,\tau)}{\psi_{\tau}(z,\bar z,\tau)}
\end{equation}
This representation is a special case of (\ref{eq: how to indicate isothermic}), 
which implies 
\begin{coro}               			\label{Corollary of isothermicness}
Bonnet surfaces are isothermic away from umbilic points.
\end{coro}

Representation (\ref{isothermicQ}) implies that $\psi_\tau$ is a harmonic function on 
${\cal R}\setminus {\cal U}$. Moreover, from a more detailed analysis of the  local behaviour
at an umbilic point one can deduce \cite{BE2} that the function $\psi_\tau$ can be
extended to a nowhere vanishing harmonic function on $\cal R$
$$
\psi_\tau:{\cal R}\to \BbbR_*.
$$

With so defined $\psi_{\tau}$, the identity
\begin{equation}\label{varphi equal something}
  \varphi = {\rm i}\,\psi_{\tau}\, Q
\end{equation}
holds on all ${\cal R}$.

\begin{theo}						\label{L:U is a subset of V}
Let ${\cal F}$ be a Bonnet surface.
Then
\begin{itemize}
\item[$(i.)$]{Umbilic points are critical points of the mean curvature function ${\cal U}\subset{\cal V}$.}
\item[$(ii.)$]{The set ${\cal V}$ of critical points of the mean curvature function 
is discrete in ${\cal R}$.}
\end{itemize}
\end{theo} 
{\it Proof}: From the Codazzi equations it follows that $d H=0$ if and only if $Q_{\bar z}=0$.
Differentiating (\ref{varphi equal something})  with respect to $\bar z$ one obtains
\begin{equation*}
        {\rm i} \,\psi_{\tau \bar z}\,Q + {\rm i}\, Q_{\bar z}\equiv 0.
\end{equation*}
Thus $Q=0$ implies $Q_{\bar z}=0$.
To show that ${\cal V}$ is a discrete subset of 
${\cal R}$, we use
\begin{equation*}
        Q_{\bar z}=0 \Longleftrightarrow
        {\rm i}\,\psi_{\tau \bar z}\,Q = 0 \Longleftrightarrow
        {\rm i}\,\psi_{\tau z}\,Q =0 \Longleftrightarrow
        \dfrac{\psi_{\tau z}\, Q_{\tau}}{\psi_{\tau}} =0
\end{equation*}
where we use that $\psi_{\tau}$ is a never vanishing harmonic
function on ${\cal R}$. Since $\psi_{\tau z}\, Q_{\tau}\,d z^3$ is a
cubic holomorphic differential, its zeros (which comprise the set ${\cal V}$) are discrete.

\subsection{Local Theory away from Critical Points}			\label{s4.2}

In this section we develop  local theory of Bonnet surfaces
in $\BbbR^3$ away from possible critical points of the mean curvature function 
$F:(-\epsilon,\epsilon)\times{\cal R}\setminus{\cal V}\to \BbbR^3$. 
The preimages of holomorphic local charts $z:U\subset{\cal R}\setminus{\cal V}\to \BbbC$
are always assumed to be simply connected. 
 
 The following ``stationary'' characterisation of Bonnet surfaces is classical 
\cite{Gra}.
\begin{theo}				\label{Lemma: 1/Q harm. iff Bonnet}
An umbilic free surface ${\cal F}$ is a Bonnet surface if and only if

\noindent (i.) ${\cal F}$ is isothermic.

\noindent(ii.) $1/Q$ is harmonic, i.e.
\begin{equation}\label{1/Qzzb0}
     \left(\dfrac{1}{Q(z,\bar z)}\right)_{z \bar z} =0,
\end{equation}
where $z$ is an isothermic coordinate and $Q\,d z^2$ is the Hopf differential.
\end{theo}
{\it Proof}: ``$\Rightarrow$'':
follows from Corollary \ref{Corollary of isothermicness}, (\ref{isothermicQ}),
and harmonicity of $\psi_\tau$.
\newline
``$\Leftarrow$'': Let $z$ be an isothermic coordinate and 
$e^{u(z,\bar z)}, H(z,\bar z), Q(z,\bar z)$ be a solution to the 
Gauss-Codazzi equations (\ref{complex version of Gauss-Codazzi}) with $Q$ satisfying 
(\ref{1/Qzzb0}).
Locally there exists  a holomorphic function
$h(z)$ such that $Q(z,\bar z) = {1/(h(z) +\bar h(\bar z))}$.
Define $Q(z,\bar z,\tau)$ via 
\begin{equation}\label{Q for the family in R-V}
        Q(z,\overline{z},\tau) = 
        \left(
        \dfrac{1 - {\rm i}\,T\,\overline{h(z)}}{1 + {\rm i}\,T\,h(z)}
        \right)\,\dfrac{1}{h(z)+\overline{h(z)}},
\end{equation}
where $T$ is a deformation parameter equivalent to $\tau$. One can easily check that
 $\vert Q(z,\bar z,\tau)\vert =
\vert Q(z,\bar z)\vert $ and $Q_{\bar z}(z,\bar z,\tau) = Q_{\bar z}(z,\bar z)$ holds. 
Thus $e^{u(z,\bar z)}$, $H(z,\bar z)$, $Q(z,\bar z,\tau)$ is also a solution of 
(\ref{complex version of Gauss-Codazzi}) for all 
$T$.  The surfaces corresponding to different $T$ are isometric and have the
same mean curvature. They form a Bonnet family.

It is easy to see that the Codazzi equations imply
\begin{equation*}
        h'(z) \, H_z(z,\bar z) =  \overline{h'(z)}\, H_{\bar z}(z,\bar z).
\end{equation*}
Introducing locally the new conformal coordinate
\begin{equation}\label{define w at non-umbilics}
        w = \int\dfrac{1}{h'(z)} \, d z
\end{equation}
one finds that the mean curvature function depends on 
\begin{equation*}
     t = w + \bar{w}
\end{equation*}
only. This finally leads to the fact that the Gauss-Codazzi equations can be reduced to an
ordinary differential equation, which is derived below. 

One can directly check that $Q$ satisfies 
\begin{equation*}
        Q_{\bar w}(w,\bar w,T) = \bar Q_w(w,\bar w,T) = 
        -\vert Q(w,\bar w,T)\vert^2.
\end{equation*}
Inserting 
\begin{equation}
  e^{u} = \dfrac{2\,Q_{\bar w}}{H_w} = - \dfrac{2 \,\vert Q\vert^2}{H'}
\end{equation}
into the Gauss equation one  obtains 
\begin{equation}				\label{Hazzidakis for non-umbilics I}
        \left(\dfrac{H''(t)}{H'(t)}\right)' -  H'(t) =
        \vert Q\vert^2\,\left(2 -\dfrac{H^2(t)}{H'(t)}\right),
        \quad ' =\dfrac{d}{d\,t}.
\end{equation}
For a general holomorphic function $h(w)$, equation
(\ref{Hazzidakis for non-umbilics I}) does not possess a solution depending on
$t$ only.

The identity $\vert Q\vert^2_w=\vert Q\vert^2_{\bar w}$
implies 
$$
        \left(\dfrac{h''(w)}{h'(w)} - \overline{\dfrac{h''(w)}{h'(w)}}
        \right)\,
        \left(h(w) +\overline{ h(w)}\right) = 2\,(h'(w) - \overline{h'(w)}).
$$
This equation can be solved explicitly. Up to appropriate normalizations $h(w)$ is one of the
following five forms
$$
\begin{matrix}
        h_1(w) & = & w, &
        h_2(w) & = &  e^{4\,{\rm i}\,w}, &
        h_3(w) & = & -\dfrac{1}{w},\cr 
        h_4(w) & = &  \tanh(2 \,w),& 
        h_5(w) & = &  \tan(2\,w). & & & 
\end{matrix}
$$

Finally one arrives at the following classical result of E. Cartan \cite{Car} (see \cite{BE1,
BE2} for detail).
\begin{theo}						\label{theoremCartan}
Away from umbilic points  there are 3 types of Bonnet families, which are
characterized by the modulus of the Hopf differential:
\smallskip 
    \newline{\bf Type A :} 
    $\vert Q^A(w,\bar w,T)\vert^2 = \dfrac{4}{ \sin^2( 2 t )}$.
    \newline{\bf Type B :} 
    $\vert Q^B(w,\bar w,T)\vert^2 = \dfrac{4}{\sinh^2(2 t)}$.
    \newline{\bf Type C :} 
    $\vert Q^C(w,\bar w,T)\vert^2 = \dfrac{1}{t^2}$ .\\[0.3cm]
The families are given by the surfaces with the fundamental forms presented in Table
\ref{tab1}.
\setlength{\extrarowheight}{4pt}
\begin{table}[thb]
\begin{center}
\begin{tabular}{| >{\footnotesize}c|| >{\footnotesize\parbox[c][1cm][c]{4cm}}c | 
>{\footnotesize}c| >{\footnotesize}c|}\hline
        Type  &    $Q(w,\bar{w})$  &  $H$ solves  & $e^u$  
        \\[0.1cm] \hline\hline 
        $A_1$ & 
                $-2\,
                \left(\dfrac{\sin ( 2 \bar w)}{\sin ( 2 w) }\right) 
                                    \dfrac{ 1}{\sin(2\,t)}$ 
              & 
              &     \\ \cline{1-2}
        $A_2$ &    $2\,\left(\dfrac{ 
                \cos ( 2 \bar w)}{\cos ( 2 w)}\right) 
                                    \dfrac{1}{\sin(2\,t)}  $ 
              & \raisebox{2.5ex}[-2.5ex]{$\left(\left( 
                \dfrac{H''}{H'}\right)' - H'\right) 
                \dfrac{\sin^2(2 t)}{4} 
                = 2 - \dfrac{ H^2}{H'}$}
              &\raisebox{2.5ex}[-2.5ex]{$-\dfrac{8}{\sin^2(2 t) H'}$}\\ \hline
    B     &    $  - 2\,
                \left(\dfrac{ \sinh ( 2 \bar w)}{\sinh (2 w)}\right)\,
                                    \dfrac{ 1}{\sinh( 2\,t)} 
               $ 
              &$\left(\left( 
                \dfrac{ H''}{H'}\right)' - H'\right) 
                \dfrac{ \sinh^2(2 t)}{4} 
                = 2 - \dfrac{H^2}{H'}$ 
              &$-\dfrac{8}{\sinh^2(2 t) H'}$ 
    \\ \hline 
    $C$ &    $-\dfrac{\bar w}{w} 
                        \dfrac{1}{t}$ 
                & $\left(\left( 
                \dfrac{H''}{H'}\right)' - H'\right) t^2 
                = 2 - \dfrac{H^2}{H'}$
                & $-\dfrac{2}{t^2 H'} $ 
                \\ \hline
\end{tabular}\caption[tab1]{\label{tab1}Table of fundamental functions}
\end{center}
\end{table}
\newline
Here $H(t)$ is any smooth solution with $H'<0$ of the corresponding ordinary differential 
equation in Table \ref{tab1}.\newline
The corresponding one-parameter families of isometries are intrinsic
isometries of the surface described by imaginary translations
of the coordinate $w$
\begin{equation}\label{isometryaction}
        w\to w + {\rm i}\, \rho(T).
\end{equation}
The surfaces of type $A_1$ and $A_2$ are isometric with the same mean
curvature function.
\end{theo}

\begin{coro}						\label{Cor p 25}
Bonnet surfaces are real analytic.
\end{coro}

The equations for the mean curvature function in Table \ref{tab1}
were first derived by N.Hazzi\-da\-kis in \cite{Haz}.
We call them Hazzidakis equations.

A Bonnet surface in $\BbbR^3$ is said to be {\em of type} {\bf A},
{\bf B}, or {\bf C}, respectively, if away from critical points 
of the mean curvature function it is of the corresponding type.
One can show that  
a Bonnet surface is  exactly of one of the types {\bf A}, {\bf B}, or {\bf C}. 

\subsection{Local Theory at Critical Points}			\label{s4.3}

In this section we derive a differential equation describing 
Bonnet surfaces near (isolated) critical points.

The identity (\ref{varphi equal something}) at an umbilic point implies
the following local representation of the Hopf differential.
\begin{prop}
Let  $P\in{\cal U}\subset {\cal R}$ be an umbilic point of a Bonnet surface.
Then there exists a neighbourhood $U$ of $P$, a local conformal chart
$z:U\to \BbbC$ with $z(P)=0$, a holomorphic 
non-vanishing function $\varphi:U\to \BbbC_*$, and an integer
$J>0$ such that the Hopf differential on $U$  is 
\begin{equation}					\label{Q at umbilic II}
        Q(z,\bar z,\tau) = -{\rm i}\,z^J\, \left( 
        \dfrac{\varphi(z,\tau)}{\psi_{\tau}(z,\bar z,\tau)}
        \right).
\end{equation}
\end{prop}

The function $\psi_{\tau\, z}$ is holomorphic on $U$ and therefore can be represented as
$\psi_{\tau z} = z^M\,\theta(z),\, \theta(0)\not=0$ with some $M\ge 0$. Analysing the 
Codazzi equations 
\begin{equation}\label{Codazzi at Umbilic}
        \bar Q_z = \bar z^J\,z^M\,\theta\,
        \overline{
        \left(
        \dfrac{{\rm i}\,\varphi}{\psi_{\tau}^2}
        \right)} = \dfrac{1}{2}\,H_{\bar z}\,e^u,\quad
        Q_{\bar z} = z^J\,\bar z^M\,\bar\theta\,        
        \left(
        \dfrac{{\rm i}\,\varphi}{\psi_{\tau}^2}
        \right) = \dfrac{1}{2}\,H_z\,e^u.
\end{equation}
one obtains 
\begin{equation}					\label{eq:define M}
        \psi_{\tau z} = z^{J+1}\,\theta(z), \qquad\theta(0)\not=0.
\end{equation}
Non-umbilic critical points correspond to $J=0$.
We call $-J$ with $J$ given by (\ref{Q at umbilic II}, \ref{eq:define M})
the {\em index of the critical point}. 

The derivation of an ordinary differential equation at a critical point is similar to the
one of section \ref{s4.2}. 
The Codazzi equations (\ref{Codazzi at Umbilic}) imply
\begin{equation}\label{Codazzi for z, II}
        \left(\dfrac{\theta}{{\rm i}\,\varphi}
        \right)\,z \,H_z =
        \overline{\left(
        \dfrac{\theta}{{\rm i}\,\varphi}        
        \right)}\,\bar z\,H_{\bar z}\in \BbbR.     
\end{equation}
Introducing a new conformal parameter $w(z)$ by 
\begin{equation}					\label{define w at umbilics}
  \dfrac{z}{f}\; \dfrac{\partial}{\partial \,z} = 
  \dfrac{w}{f(0)}\;\dfrac{\partial}{\partial \,w},
\end{equation}
with $w(0)=0$ and $f(z) =-i\,\varphi(z)/\theta(z)$, one arrives at
\begin{equation*}
  w\, H_w = \bar w\, H_{\bar w}.
\end{equation*}
The last identity shows that the mean curvature function is a function of 
$s=\vert w\vert^2$ only. Proceeding further as in Section \ref{s4.2}, after apropriate
normalizations (see \cite{BE2} for detail) one obtains the following

\begin{theo}\label{Theorem Umbilics M3(c)}
Let ${\cal F}$ be a Bonnet surface in $\BbbR^3$ with an
(isolated) critical point of index $-J$. Then there exists a local conformal
chart $w$  at the critical point $w=0$  such that the 
coefficients of its fundamental forms are given by
\begin{equation}\label{FundFunctatInfinity}
\begin{matrix}
        Q(w,\bar w)\, d w^2 & = & (J+2)\,\left(
        \dfrac{1 - \bar w^{J+2}}{ 1 - w^{J+2}}\right)\,
        \dfrac{w^J}{ 1 - s^{J+2}}\, d w^2,
        \vspace{2mm}\cr
        e^{u(w,\bar w) }\,\vert d w\vert^2 & = &
        -2 \,(J+2)^2\,\dfrac{s^J}{ ( 1- s^{J+2})^2\,H'(s)}
        \,\vert d w\vert^2,
\end{matrix}
\qquad s = \vert w\vert^2\;,
\end{equation}
and $H(w,\bar w) = H(\vert w\vert^2)$ is a solution of
\begin{equation}\label{eq: H at umbilics}
        \left(\dfrac{s H''(s)}{H'(s)}\right)' - H'(s) = 
        \dfrac{(J+2)^2\,s^{J+1}}{(1- s^{J+2})^2}\,\left(
        2 - \dfrac{H^2(s)}{s\,H'(s)}\right),
      \quad
      ' =\dfrac{d}{d\,s },
\end{equation}
with the asymptotics
\begin{equation}\label{eq:asymptotics for H at s eq to 0}
  H(s) = H(0) + s^{J+1}\,B(s)
\end{equation}
at $s=0$, where $B(s)$ is a non-vanishing, smooth 
function.\newline
Conversely, any solution of (\ref{eq: H at umbilics}) with the asymptotics
(\ref{eq:asymptotics for H at s eq to 0}) at $s=0$ via 
(\ref{FundFunctatInfinity}) 
determines a Bonnet surface with critical point of index $-J$. The Bonnet family is given 
by the intrinisic isometry 
\begin{equation}\label{intrinisc isometry for the coord}
     w\to e^{{\rm i}\,\alpha}\,w,\quad\alpha\in\BbbR.
\end{equation}
\end{theo} 

The existence of Bonnet surfaces in $\BbbR^3$ with critical points of the
mean curvature function will be proven in Section \ref{s4.5}.

Finally note that Bonnet surfaces with critical points are necessarily of type {\bf B}.
The relation between the corresponding coordinates $w$ of this section and $w_{\bf
B}=w$ of Section \ref{s4.2} is given by
\begin{equation}\label{wBV to wB}     
        w = e^{-\frac{4}{J+2}\,w_{\bf B}}.
\end{equation}
We call a Bonnet surface with an isolated critical point a Bonnet surface of type 
${\bf B}_{\cal V}$.

\subsection{Bonnet Surfaces via Painlev\'{e} Transcendents}		\label{s4.4}

In this section we give an explicit description of Bonnet surfaces in terms of solutions 
of Painlev\'{e} equations, which are certain ordinary differential equations of the 
second order
\begin{equation}							\label{2dODE}
y''= R(y', y, x).
\end{equation} 

Solutions of the Painlev\'{e} equations -- Painlev\'{e} transcendents -- are treated nowadays 
as non-linear 
special functions. The theory of these special functions is rather well developed
(see for example \cite{IN, IKSY, Its}). 
In Section \ref{s4.5} we use elements of this theory for global classification of Bonnet
surfaces. Of special importance for us is the {\em Painlev\'{e} property} of these equations.
Recall that a differential equation (\ref{2dODE}) is said to possess the Painlev\'{e} property 
if it is free of movable branch and essential singular points, i.e. the only singularities of 
the solutions which change their position if one varies the initial data are poles.

Regarding examples of the Bonnet surfaces of type {\bf B}, let us expain how the quaternionic 
representation of Bonnet surfaces naturally leads to the above mentioned remarkable 
connection. Note that the Hazzidakis equations of types {\bf A} and {\bf B} are analytically
equivalent $H_{\bf A}(it)=-i H_{\bf B}(t)$.

Since the mean curvature and the metric of a Bonnet surfaces of type {\bf B} depend on 
$t=w+\bar w$ only, it is natural to change the parametrization once more using $t$ 
(or a function of $t$) as one of the variables.
After some computations one arrives to the idea of considering
\begin{equation}\label{def: new variables to reach Fuchsian system}
        x= e^{-4\,(w+\bar w)},
        \qquad \lambda = e^{-4 \,w}
\end{equation}
as new variables. In these coordinates the Hazzidakis equation is 
\begin{equation}\label{B:Hazzidakis for non-umbilics II}
        4\,\left(x\,\dfrac{{\cal H}''(x)}{{\cal  H}'(x)}\right)' +
        {\cal H}'(x) = \dfrac{4}{(x-1)^2}\,
        \left(
        2 + \dfrac{{\cal H}^2(x)}{4\, x\,{\cal H}'(x)}\right),
\end{equation}
where ${\cal H}(x) \equiv H(t)$.
The frame equations 
(\ref{R32.2.16}, \ref{R32.2.18}) 
in parametrization (\ref{def: new variables to reach Fuchsian system}) are now
\begin{eqnarray}
        \dfrac{\partial \Phi}{\partial \lambda}\,\Phi^{-1} &=&
        \dfrac{B_{0}(x)}{\lambda} +
        \dfrac{B_{1}(x)}{\lambda-1} +
        \dfrac{B_{x}(x)}{\lambda-x},						\label{lambda-eq}\\
        \dfrac{\partial \Phi}{\partial x}\,\Phi^{-1} &=&
        - \dfrac{B_{x}(x)}{\lambda-x} + C(x), 					\label{x-eq}
\end{eqnarray}
where all the coefficients of the matrices are given by some explicit formulas through 
${\cal H}(x), {\cal H}'(x)$ and ${\cal H}''(x)$. 
Ignore for the moment the complicated formulas for the coefficients in the
matrices (\ref{lambda-eq}). What is more important is the special dependence of (\ref{lambda-eq}) 
on $\lambda$.
Equation (\ref{lambda-eq}) is a 2 by 2 matrix dimensional Fuchsian system with four regular 
singularities (at $\lambda=0, 1, x, \infty$). Equation (\ref{x-eq}) describes 
deformations preserving the 
monodromy group of the system. It is well known (see for example \cite{JM}) that such isomonodromic 
deformations are the issue of the Painlev\'{e} VI equation.
In particular, all the coefficients of the matrices can be expressed in terms of solutions of this 
equation -- Painlev\'e transcendents.
A rather involved calculation (see \cite{BE1}) of the corresponding gauge transformation identifying the
corresponding descriptions yields finally the following 
\begin{theo}\label{B:Mapping y to H and inverse}

Equation (\ref{B:Hazzidakis for non-umbilics II}) possesses the first integral
\begin{equation}					\label{B:1.Integral for H}
        x^2 \,\left(\dfrac{{\cal H}''(x)}{{\cal H}'(x)} +\dfrac{2}{x-1}\right)^2 
        + \dfrac{x\,{\cal H}'(x)}{2}
        +\dfrac{{\cal H}^2(x)}{2\,(x-1)^2\,{\cal H}'(x)} + 
        \dfrac{{\cal H}(x)\,(x +1)}{2\,(x-1)} 
        = \theta^2.
\end{equation}

Let ${\cal H}(x)$ be a solution of (\ref{B:Hazzidakis for non-umbilics II})
with ${\cal H}'(x)\not=0$
and $\theta$ be a fixed root of (\ref{B:1.Integral for H}). Then 
the function $y(x)$ defined by
\begin{equation}\label{def y in H for Bonnet B}
        y(x) \equiv - \dfrac{2}{{\cal H}'(x)}\,\left(
        \dfrac{x\,(x -1)\,{\cal H}''(x)+ 
        (\theta - x\,(\theta -2))\,{\cal H}'(x)}{ 
        {\cal H}(x)+ (x - 1)\,{\cal H}'(x)}
        \right)^2
\end{equation}
solves the Painlev\'e VI equation of the following form
\begin{equation}\label{PVI for Bonnet B type}
\begin{matrix}
\dfrac{d^2\, y}{d\,x^2} & =  & 
        \dfrac{1}{2}
        \left( 
        \dfrac{1}{y} + 
        \dfrac{1}{y - 1} +
        \dfrac{1}{y - x}
        \right)\,{y'}^2 -
        \left(
        \dfrac{1}{x}+
        \dfrac{1}{ x-1} +
        \dfrac{1}{ y - x}
        \right)\,y' +\cr
        & & 
        \dfrac{y\,(y-1)\,(y-x)}{2\,x^2\,(x-1)^2}\,
        \left(
        \theta^2\,\dfrac{(x-1)}{(y-1)^2} 
        -\theta\,(\theta +2)\,\dfrac{x\,(x-1)}{(y-x)^2} 
        \right).
\end{matrix}
\end{equation}
Conversely, for any solution of the Painlev\'e VI (\ref{PVI for Bonnet B type}) which 
is not of the form 
 $       y(x) = {\rm const}\cdot\,x^{-\theta}$, 
the function 
\begin{equation}\label{def H in y for Bonnet B}
        {\cal H}(x) \equiv 
        -2 \,\dfrac{(x -1)\,(\theta^2\,y(x)^2 -x^2{y'}^2(x))}
        {y(x)\,(y(x) -1)(y(x) -x)} 
\end{equation}
is a solution of (\ref{B:Hazzidakis for non-umbilics II}) with first
integral $\pm \theta$ (\ref{B:1.Integral for H}).

The mappings (\ref{def y in H for Bonnet B}) and (\ref{def H in y for Bonnet B}) with the same 
$\theta$ are inverse one to another.
\end{theo}

In a similar way (see \cite{BE1} for the corresponding formulas), Bonnet surfaces of type {\bf C} 
are described through solutions of the Painlev\'e V equation.

\begin{coro}					\label{Cor: Painleve Property}
The Hazzidakis equations for Bonnet surfaces of 
all types possess the Painlev\'e property.
\end{coro}

\subsection{Global Classification of Bonnet Surfaces}			\label{s4.5}

In this section we classify all {\em maximal} or {\em global} Bonnet immersions 
$F:{\cal R}\to\BbbR^3$.
These surfaces are characterized by the following natural ``maximality'' property:
given an immersed Bonnet surface ${\cal F}\subset \BbbR^3$ there exists a
global Bonnet immersion containing ${\cal F}$, i.e. 
\begin{equation*}
  {\cal F} = F(U), \quad U\subset{\cal R}.
\end{equation*}

Let us first prove the existence of critical points.
\begin{theo}							\label{Corollary on BV}
For arbitrary $H(0)\in\BbbR$, $H_0<0$ there exists a real analytic 
Bonnet surface of type ${\bf B}_{\cal V}$
with a critical point of index $-J$ and with the mean curvature function and the metric at
the critical point given by
\begin{equation*}
  H(0)\mbox{ and } -2\,\dfrac{(J+2)^2}{(J+1)\,H_0}\,d w\,d \bar w.
\end{equation*}
\end{theo}
{\it Proof.} Substituting the ansatz 
$$
H(s) = H(0) + s^{J+1} \,{\sum_{i=0}^{\infty}}\,H_i\, s^i,\quad
        H_0\not=0
$$
into equation (\ref{eq: H at umbilics}) and using the corresponding recurrence formulas 
for $H_i$, one
can prove that the series converges absolutely in a neighbourhood of $s=0$. Thus,
for arbitrary $H(0)\in\BbbR$, $H_0\in\BbbR_*$ there exists a real analytic  solution  of 
equation (\ref{eq: H at umbilics}) at $s=0$  with the asymptotic
\begin{equation}					\label{formal serie for H}
  H(s) = H(0) + s^{J+1} H_0 + O(s^{J+2}).
\end{equation}
The claim follows now from Theorem \ref{Theorem Umbilics M3(c)}.

Let the local coordinate $w$ be defined by (\ref{define w at non-umbilics}) for the 
Bonnet surfaces of type {\bf A}, {\bf B}, and {\bf C} 
and by (\ref{define w at umbilics}) for the Bonnet surfaces of type
${\bf B}_{\cal V}$. We denote by ${\bf U}$  the largest connected domain in 
the $w$-plane, for which $Q$ (see Table {\bf\ref{tab1}} and 
(\ref{FundFunctatInfinity}) respectively) is bounded. Furthermore,
let ${\cal D}=\{t=w +\bar w \,\vert w\in{\bf U}\}$ for the types {\bf A}, {\bf B},
and {\bf C}, and  ${\cal D}=\{s=\vert w\vert^2\,\vert w\in{\bf U}\}$ for the type
${\bf B}_{\cal V}$. In the following Table {\bf\ref{tab2}}, all the cases
${\bf U}$, ${\cal D}$, and the harmonic function $\psi_{\tau}$, parametrized by $w$,
are listed:
\setlength{\extrarowheight}{4pt}
\begin{table}[h]
\begin{center}
\begin{tabular}{| c || c | c | c |}
\hline
Type  &  ${\bf U}$  & ${\cal D}$ &  $\psi_{\tau}$ \\[0.1cm] \hline\hline
        ${\bf A}_1$,  ${\bf A}_2$
        & 
        $\{w\in\Bbb{C}\;\vert\; 0<{\rm Re}(w)<\dfrac{\pi}{4}\}$ &
        $(0,\dfrac{\pi}{2}) $ &
        $2\, {\rm Re}(\tan(2 w))$\\\hline
        ${\bf B}$ & 
        $\{w\in\Bbb{C}\;\vert\;{\rm Re}(w)>0\}$ &
        $ (0,\infty)$  &
        $2\, {\rm Re}(\tanh(2 w) )$\\\hline 
        ${\bf B}_{\cal V}$ & 
        $\{w\in\Bbb{C}\;\vert\;\;\vert w\vert <1\}$ &
        $[0,1)$ & 
        $2\,{\rm Re}
        \left( \dfrac{1 + w^{J+2}}{1-w^{J+2}}\right)$\\\hline
        ${\bf C}$ &
        $\{w\in\Bbb{C}\;\vert\;{\rm Re}(w)>0\} $ &
        $ (0,\infty)$ &
        $2\,{\rm Re}(w)$\\\hline
\end{tabular}\caption[tab2]{\label{tab2}Global description of Bonnet surfaces}
\end{center}
\end{table}


\begin{prop}							\label{theorem A}
Let $H$ be a solution of one of the Hazzidakis equations of types {\bf A}, {\bf B},
{\bf C}, or ${\bf B}_{\cal V}$ (see Table {\bf\ref{tab1}} and (\ref{eq: H at umbilics}))
with $H'(t)<0$ at some point $t\in{\cal D}$ (the corresponding domains ${\cal D}$
are listed in Table {\bf\ref{tab2}}). Then $H$ is real analytic on ${\cal D}$.
\end{prop}

The proof of this proposition is based on using the Painlev\'e property 
(see Corollary \ref{Cor: Painleve Property}) together with the following 

\begin{lemm}						\label{Theorem H exists globally}
Let ${\cal D}\ni t$ be an open interval with smooth
positive-valued functions $f,\,g,\,\vert Q\vert^2:{\cal D}\to\BbbR_+$ on it. Let 
$H=H(t)$ be a real-valued solution of the generalized Hazzidakis
equation
\begin{equation}\label{the most general Hazzi}
        \left(f(t)\,\dfrac{H''(t)}{H'(t)}\right)' - H'(t) =
        \vert Q\vert^2\,\left(
        2 - \dfrac{H^2(t)}{g(t)\,H'(t)}\right),
\end{equation}
smooth on ${\cal D}\setminus{\cal P}$, where ${\cal P}$ is a discrete set
of poles of $H(t)$.
\newline
If $H'(t_0)<0$ at some $t_0\in{\cal D}$, then $H$ is smooth everywhere on 
${\cal D}$, i.e. ${\cal P}=\emptyset$, with $H'(t)<0$ for all $t\in{\cal D}$.
\end{lemm}

To prove smoothness of $H(t)$ one shows that all poles of $H(t)$ are necessarily 
simple with negative residues and that in addition $H'(t)$ never changes its sign. These two properties
contradict one another.

It is not difficult to show that immersions $F:{\bf U}\to\BbbR^3$ of Bonnet surfaces of type {\bf A}, {\bf B},
{\bf C}, or ${\bf B}_{\cal V}$ given in Tables \ref{tab1} and \ref{tab2} are maximal. 
The function $\psi_{\tau}$ defined in Table \ref{tab2} is a non-vanishing
function on $\bf U$. The continuity of $\psi_{\tau}$ yields that this function can not be extended beyond
$\bf U$. 

Finally using the arguments of Section \ref{s4.2} it is easy to show that the classified global Bonnet 
surfaces are all different, i.e. for 
two global Bonnet immersions $F_i:{\cal R}_i\to\BbbR^3$, $i=1,2$ there exist
no open sets $U_i\subset {\cal R}_i$, $i=1,2$ on which the surfaces 
coincide.

\begin{theo}
Any Bonnet surface in $\BbbR^3$ can be conformally parametrized 
$F:{\cal R}\to \BbbR^3$  by a corresponding global Bonnet immersion
$F:{\bf U}\to\BbbR^3$, ${\cal R}\subset {\bf U}$. The latter are of one of the types {\bf A},
{\bf B}, {\bf C}, ${\bf B}_{\cal V}$. The corresponding domains $\bf U$ are listed
in Table \ref{tab2}.\newline
Given $t_0\in{\cal D}$ (see Table \ref{tab2}) and arbitrary $H(t_0)$,
$H'(t_0)<0$, $H''(t_0)$ there exists a unique solution $H(t)$
of the Hazzidakis equation of type {\bf A}, {\bf B}, {\bf C} (see Table \ref{tab1}),
real analytic on {\cal D}. This function determines the fundamental forms (Table \ref{tab1}) of the
corresponding global Bonnet immersions of the type {\bf A}, {\bf B}, {\bf C}.
Given $H(0)$, $H_0<0$ there exists a  unique solution $H(s)$ of the Hazzidakis equation 
(\ref{eq: H at umbilics}), real analytic on ${\cal D}$, with the asymptotics 
(\ref{formal serie for H}). It determines by (\ref{FundFunctatInfinity}) the fundamental forms of the 
corresponding global Bonnet surface of type ${\bf B}_{\cal V}$.
\end{theo}

\subsection{Examples of Bonnet Surfaces}				\label{s4.6}


Let us present some figures of Bonnet surfaces of types {\bf B} and ${\bf B}_{\cal V}$.
Tubes correspond to parameter lines $t={\rm const}$, i.e.
to the trajectories of the isometric flow preserving the mean curvature function.
Both the mean curvature function and the metric are preserved along these lines.
The last fact can be clearly observed - the strips bounded by two consequent
parameter lines $t=t_1$ and $t=t_2$ are of constant width.
The isometry is intrinsic, i.e. is not induced by a Euclidean motion of the ambient $\BbbR^3$.
The immersion domain {\bf U} of Bonnet surfaces
of type {\bf B} is naturally split into fundamental domains
\begin{equation*}
  {\bf U}_n =\{w\in\Bbb{C}\,\vert \, (n-1)\dfrac{\pi}{2}< {\rm Im}(w)< n\,\dfrac{\pi}{2}\}.
\end{equation*}
Indeed the fundamental forms (see Table {\bf\ref{tab1}}) are invariant with respect to the shift
\begin{equation*}
  w\to w + {\rm i}\,\dfrac{\pi}{2},
\end{equation*}
and thus immersed ${\bf U}_n$ with different $n$'s are congruent in $\BbbR^3$. A Bonnet surface
comprised of three fundamental domains is shown in Figure 7. For an appropriate
choice of parameters, several copies of the fundamental domain can close up and thus
comprise a closed surface with a critical point. Figure 8 shows such a case
with three fundamental domains. It is worth mentioning that it was this figure which 
lead us to suggest the existence of Bonnet surfaces with critical points.

Figures 9 and 10 present another example of type 
${\bf B}_{\cal V}$. The surface in Figure 9 is an immersed disk. 
The index of the critical point is $J=6$. 
A more detailed view of one of the cusps of this surface is
shown in Figure 10. As before,
the tubed curves are integral curves of the isometry field. The surface is probably embedded.

The figures of this section are produced by Ulrich Eitner. Further examples can be found in 
\cite{BE1, BE2}.



\newpage
\begin{thebibliography}{WWW}
\bibitem[Abr]{Abr}
Abresch, U.: Constant mean curvature tori in terms of elliptic functions.
J. reine angew. Math. {\bf 394} (1987) 169-192
\bibitem[BBEIM]{BBEIM}
Belokolos, E.D., Bobenko, A.I., Enolskii, V.Z., Its, A.R., Matveev, V.B.:
Algebro-Geometric Approach to Nonlinear Integrable Equations, Springer, Berlin
(1994)
\bibitem[Bia1]{Bianchi}
Bianchi,L.: Vorlesung \"uber Differentialgeometrie,
Leipzig, Berlin, p. 454 (1910)
\bibitem[Bia2]{BianchiBonnet}
Bianchi, L.: Sulle superficie a linee di curvature isotherme.
Rend. Acc. Naz. dei Lincei, (5) {\bf  12} (1903), 511-520
\bibitem[Bia3]{BianchiIsoth}
Bianchi, L.: Ricerche sulle superficie isotherme and sulla deformatione delle
quadriche. Ammali di Math. Ser. III {\bf 11} (1905) 93-157
\bibitem[Bob1]{BoCMC}
Bobenko, A.I.: Constant mean curvature surfaces and integrable equations.
Russ. Math. Surveys {\bf 46} (4) (1991), 1 - 45
\bibitem[Bob2]{Bo2x2}
Bobenko, A.I.:
Surfaces in Terms of 2 by 2 Matrices. Old and New Integrable Cases.
In: Fordy, A.P., Wood, J.(eds.)
Harmonic Maps and Integrable Systems, Braunschweig Wiesbaden: Vieweg (1994)
83-127
\bibitem[BE1]{BE1}
Bobenko, A.I., Eitner, U.:
Bonnet surfaces and Painlev\'e equations.
J. reine angew. Math. {\bf 499} (1998), 47-79
\bibitem[BE2]{BE2}
Bobenko, A.I., Eitner, U.:
Painlev\'e equations in differential geometry of surfaces (in preparation).
\bibitem[BP]{BP}
Bobenko, A., Pinkall, U.: Discrete isothermic surfaces. J. reine angew. Math.
{\bf 475} (1996) 187-208
 \bibitem[Bon]{Bonnet}
Bonnet, O.: M\'emoire sur la th\'eorie des surfaces applicables,
J. Ec. Polyt. {\bf 42} (1867), 72-92
\bibitem[Car]{Car}
Cartan, E.:
Sur les couples de surfaces applicables avec conversation des courbieres
principles, Bull. Sc. Math {\bf 66},  1-30 (1942)
\bibitem[Che]{Che} 
Chern, S.S., 
Deformation of Surfaces Preserving Principal Curvatures. 
In: Chavel, I., Farkas, H., (eds.), 
Differential Geometry and Complex Analysis , 
Springer 1985, 155-163 
\bibitem[Cie]{Cie}
Cieslinski, J.: The Darboux-Bianchi transformation for isothermic surfaces.
Classical results versus the soliton approach. Diff. Geom. and its Appl. {\bf 7}
(1997) 1-28
\bibitem[Dar]{Dar}
Darboux, G.: Sur les surfaces isothermiques. Comptes Rendus {\bf 122} (1899)
1299-1305, 1483-1487, 1538
\bibitem[DKN]{DKN}
Dubrovin, B.A., Krichever, I.M., Novikov, S.P.: Integrable systems I. In:
Contemprorary problems of mathematics, Fundamental Directions, Itogi nauki i Tekhniki, VINITI AN
SSSR, Moscow, vol. 4 (1985) 210-315
\bibitem[DPW]{DPW}
Dorfmeister, J., Pedit, F., Wu, H.:
Weierstrass type representation of harmonic maps into symmetric spaces.
GANG preprint III. {\bf 25} (1994) (to appear in Comm. in Anal. Geom.)
\bibitem[Eis]{Eis}
Eisenhart, L.P.: A Treatize on the Differential Geometry of Curves and Surfaces.
Allyn and Bacon, Boston (1909)
\bibitem[FPPS]{FPPS}
Ferus, D., Pedit, F., Pinkall, U., Sterling, I.:
Minimal tori in $S^4$. J. reine angew. Math. {\bf 429} (1992) 1-47 
\bibitem[Gra]{Gra}
Graustein, W.C.:
Applicability with Preservation of Both Curvatures.
Bull. Amer. Math. Soc. {\bf 30} (1924), 19-23 
\bibitem[GKS]{KGBKS} 
Grosse-Braukmann, K., Kusner, R., Sullivan, J.: Classification of embedded
constant mean curvature surfaces with genus zero and three ends. Preprint (1997) 
\bibitem[Haz]{Haz}
Hazzidakis, J.N.:
Biegung mit Erhaltung der Hauptkr\"ummungsradien.
Journal {f\"ur} Reine und Angewandte Mathematik {\bf 117} (1887), 42-56
\bibitem[Hei]{Hei}
Heil, M.: Numerical tools for the study of finite gap solutions of integrable
systems. PhD Thesis, TU Berlin (1995)
\bibitem[Hit]{Hit}
Hitchin, F.S.: Harmonic maps from a 2-torus to the 3-sphere. J. Diff. Geom. {\bf
31} (1991) 627-710
\bibitem[Hop]{Hop}
Hopf, H.: Differential geometry in the large. Lecture Notes in Mathematics {\bf
1000} Springer, Berlin (1983)
\bibitem[Its]{Its}
Its, A.R.:
The Painlev\'e Transcendents as Nonlinear Special Functions.
In: Levi D., Winternitz P. (eds.) The Painlev\'e Transcendents, Their
Asymptotics and Physical Applications.
(pp. 40-60), New York: Plenum 1992
\bibitem[IN]{IN}
Its, A.R., Novokshenov, V.Y.:
The Isomonodromic Deformation Method in the Theory of Painlev\'e Equations.
Lectures Notes in Mathematics No. 1191, Berlin, 1986
\bibitem[IKSY]{IKSY}
Iwasaki, K., Kimura, H., Shimomura, S., Yoshida, M.:
From Gauss to Painlev\'e. A Modern Theory of Special
Functions, Braunschweig, Vieweg, 1991
\bibitem[Jag]{Jag}
Jaggy, C.: On the classification of constant mean curvature tori in $\BbbR^3$.
Comment. Math. Helvetici {\bf 69} (1994) 640-658
\bibitem[JM]{JM}
Jimbo, M.T., Miwa, T.:
Monodromy Perserving Deformation of Linear Ordinary Differential
Equations with Rational Coefficients, II.
Physica {\bf 2 D}, 407-448 (1981)
\bibitem[KPP]{KPP}
Kamberov, G., Pedit, F., Pinkall, U.: Bonnet pairs and isothermic
surfaces, Duke Math. J. {\bf 92} (3), (1998) 637-644
\bibitem[Kap]{Kap}
Kapouleas, N.: Compact constant mean curvature surfaces in Euclidean three-space.
J. Diff. Geom. {\bf 33} (1991) 683-715
\bibitem[Ken]{Ken}
Kenmotsu, K.: Weierstrass formula for surfaces of prescribed mean curvature.
Math. Ann. {\bf 245} (1979) 89-99
\bibitem[Kon]{Kon}
Konopelchenko, B.G.: Induced surfaces and their integrable dynamics. Stud. Appl.
Math. {\bf 96} (1996) 9-51
\bibitem[KS1]{KuS1}
Kusner, R., Schmidt, N.: The spinor representation of minimal surfaces in space, GANG
Preprint III.27 (1993)
\bibitem[KS2]{KuS2}
Kusner, R., Schmidt, N.: The spinor representation of surfaces in space, GANG
Preprint IV.18 (1996)
\bibitem[LT]{LawTri}
Lawson, H.B: Tribuzy, R. de,
On the Mean Curvature Function for Compact Surfaces,
J. Diff. Geom., 16, pp 179- 183 (1981)
\bibitem[PP]{PP}
Pedit, F., Pinkall, U.: Quaternionic analysis on Riemann surfaces and
differential geometry, Doc. Math., Extra Volume ICM 1998, 389-400
\bibitem[Pin]{Pin}
Pinkall, U.: Regular homotopy classes of immersed surfaces. Topology {\bf 24}
(1985) 421-434
\bibitem[PS]{PS}
Pinkall, U., Sterling, I.: On the classification of constant mean curvature
tori, Ann. Math. {\bf 130}, 407 - 451 (1989)
\bibitem[Rou]{Rou}
Roussos, I.M.: Global Results on Bonnet surfaces, J. Geom., to appear.
\bibitem[Sul]{Sul}
Sullivan, D.: The spinor representation of minimal surfaces in space. Notes
(1989)
\bibitem[Tai]{Tai}
Taimanov, I.: The Weierstrass representation of closed surfaces in $\BbbR^3$.
Sfb 288 Preprint {\bf 291} (1997)
\bibitem[Wal]{Wal}
Walter, R.: Explicit examples to the H-problem of heinz Hopf. Geom. Dedicata
{\bf 23} (1987) 187-213
\bibitem[Wen]{Wen}
Wente, H.C.: Counterexample to a conjecture of H. Hopf. Pac. J. Math. {\bf 121}
(1986) 193-243
\end{thebibliography}
\end{document}